\documentclass[11pt]{article}
\usepackage{amsthm,amsmath,amssymb}
\usepackage{amsfonts}%
\usepackage{graphicx}%
\usepackage{url}%
\usepackage{color}
\usepackage{a4wide}
\usepackage{algorithmic,algorithm,float}
\nonstopmode
\renewcommand{\Re}{\mathop{\mathrm{Re}}}
\renewcommand{\Im}{\mathop{\mathrm{Im}}}

\renewcommand{\i}{\mathrm{i}}

\newcommand{\CC}{{\mathbb C}}

\newcommand{\bI}{{\bf I}}

\newcommand{\bM}{{\bf M}}
\newcommand{\bN}{{\bf N}}

\newcommand{\bs}{{\bf s}}

\newcommand{\bt}{{\bf t}}

\begin{document}

\title{Towards computing the harmonic-measure distribution function for the middle-thirds Cantor set}

\author{Christopher C. Green \& Mohamed M. S. Nasser}

\date{}
\maketitle

\vskip-0.8cm %
\centerline{Wichita State University, Department of Mathematics, Statistics \& Physics,} %
\centerline{Wichita KS 67260, USA}%
\centerline{\tt christopher.green@wichita.edu, mms.nasser@wichita.edu}

\begin{abstract}
This paper is concerned with the numerical computation of the harmonic-measure distribution function, or $h$-function for short, associated with a particular planar domain. This function describes the hitting probability of a Brownian walker released from some point with the boundary of the domain.  We use a fast and accurate boundary integral method for the numerical calculation of the $h$-functions for symmetric multiply connected slit domains with high connectivity.
In view of the fact that the middle-thirds Cantor set $\mathcal{C}$ is a special limiting case of these slit domains, the proposed method is used to approximate the $h$-function for the set $\mathcal{C}$. We also numerically analyze some asymptotic features of the calculated $h$-functions.
\end{abstract}

\begin{center}
	\begin{quotation}
	\end{quotation}
\end{center}

\begin{center}
	\begin{quotation}
	\end{quotation}
\end{center}


\section{Introduction} \label{sec:int}
The so-called $h$-function, or harmonic-measure distribution function, associated with a particular planar domain $\Omega$ in the extended complex plane $\overline{\CC}=\CC\cup\{\infty\}$ can be viewed as a signature of the principal geometric characteristics of the boundary components of the domain $\Omega$ relative to some basepoint $z_0\in\Omega$. Such functions are also probabilities that a Brownian particle reaches a boundary component of $\Omega$ within a certain distance from its point of release $z_0$. 
The properties of two-dimensional Brownian motions were investigated extensively by Kakutani~\cite{ka44} who found a deep connection between Brownian motion and harmonic functions (see also~\cite{ka45,bookBM,SnipWard16}).
Stemming from a problem proposed by Stephenson and listed in~\cite{bh89}, the $h$-functions were first introduced as objects of study by Walden \& Ward~\cite{wawa96}. Since then, the theory of $h$-functions has been successively developed in several works~\cite{bawa14,beso03,gswc,SnipWard05,SnipWard08,wawa96,wawa01,wawa07}. For a overview of the main properties of $h$-functions, the reader is referred to the survey paper~\cite{SnipWard16}. 

Given a domain $\Omega$ in the extended complex plane $\overline{\CC}$ and a fixed basepoint $z_0\in\Omega$, the $h$-function is the piecewise smooth, non-decreasing function, $h\,:\,[0,\infty)\to[0,1]$, defined by
\begin{equation}\label{eq:hr}
	h(r)= \omega(z_0,\partial\Omega\cap\overline{B(z_0,r)},\Omega)
\end{equation}
where $\omega$ denotes harmonic measure and $B(z_0,r)=\{z\in\CC\,:\,|z-z_0|<r\}$, i.e., $h(r)$ is the value of the harmonic measure of $\partial\Omega\cap\overline{B(z_0,r)}$ with respect to $\Omega$ at the point $z_0$. 
Thus $h(r)$ is equal to $u(z_0)$, where $u(z)$ is the solution to the following boundary value problem (BVP) of Dirichlet-type:
\begin{subequations}\label{eq:bdv-u}
	\begin{align}
		\label{eq:u-Lap}
		\nabla^2 u(z) &= 0 \quad \mbox{if }z\in \Omega, \\
		\label{eq:u-1}
		u(z)&= 1 \quad \mbox{if }z\in \partial\Omega\cap\overline{B(z_0,r)}, \\
		\label{eq:u-0}
		u(z)&= 0 \quad \mbox{if }z\in \partial\Omega\backslash \overline{B(z_0,r)}. 
	\end{align}
\end{subequations}
This is illustrated in Figure~\ref{fig:h30} for the domain $\Omega_m$ defined in~\eqref{eq:Omeg-m} below, and we refer to $\partial B(z_0,r)$ as a `capture circle' of radius $r$ and center $z_0$.

The main purpose of this work is to utilize a fast and accurate numerical method for approximating the values of $h$-functions which are close to the actual, hypothetical, $h$-function associated with the middle-thirds Cantor set $\mathcal{C}$, which is defined by
\begin{equation}\label{eq:E}
	\mathcal{C}=\bigcap_{\ell=0}^{\infty} E_\ell
\end{equation}
where the sets $E_\ell$ are defined recursively by
\begin{equation}\label{eq:Ek}
	E_\ell=\frac{1}{3}\left(E_{\ell-1}-\frac{1}{3}\right)\bigcup\frac{1}{3}\left(E_{\ell-1}+\frac{1}{3}\right), \quad \ell\ge 1,
\end{equation}
and $E_0=[-1/2,1/2]$. It is a type of fractal domain which is infinitely connected. We will demonstrate that the $h$-function for the Cantor set $\mathcal{C}$ can be approximated by computing $h$-functions for multiply connected slit domains of finite, but high, connectivity.

For a fixed $\ell=0,1,2,\ldots$, the set $E_\ell$ consists of $m = 2^\ell$ slits $I_{j}$, $j = 1, 2, \ldots, m$, numbered from left to right.  The slits have the same length $L = (1/3)^\ell$.
The center of the slit $I_{j}$ is denoted by $w_{j}$.
The domain $\Omega_m$ is the complement of the closed set $E_\ell$ with respect to the extended complex plane $\overline{\CC}$, i.e.,
\begin{equation}\label{eq:Omeg-m}
	\Omega_m = \overline{\CC}\setminus\bigcup_{j=1}^{m} I_j, \quad m=2^\ell.
\end{equation}
As the value of $m$ is increased, the computed $h$-function for the domain $\Omega_m$ provides a successively better approximation of the actual $h$-function associated with the middle-thirds Cantor set $\mathcal{C}$.

The particular location of the basepoint $z_0$ is intrinsic to $h$-function calculations. In the current work, we consider $h$-functions for $\Omega_m$ with two different basepoint locations: one strictly to the left of all slits at $z_0=-3/2$, and another at $z_0=0$, as illustrated in Figure~\ref{fig:h30} for $m=2$. 
It should be highlighted that the method to be proposed in this paper, with modifications, can still be used to calculate $h$-functions for other basepoint locations.

\begin{figure}[ht] %
	\centerline{\hfill
		\scalebox{0.4}{\includegraphics[trim=0 0 0 0,clip]{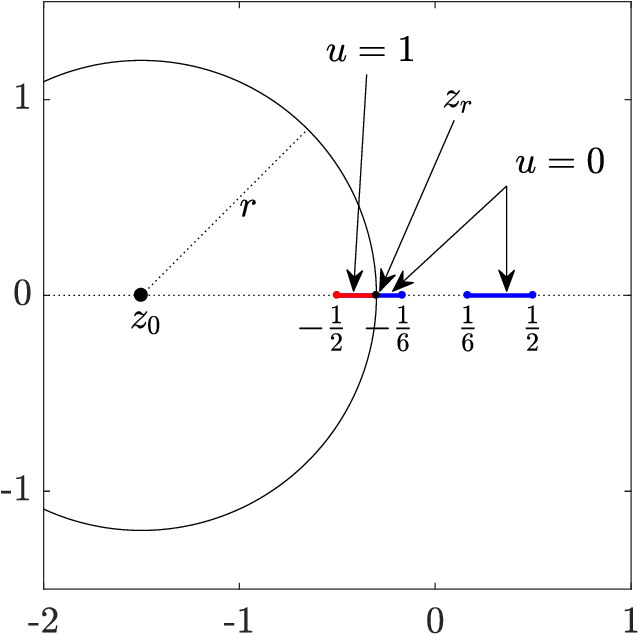}}\hfill
		\scalebox{0.4}{\includegraphics[trim=0 0 0 0,clip]{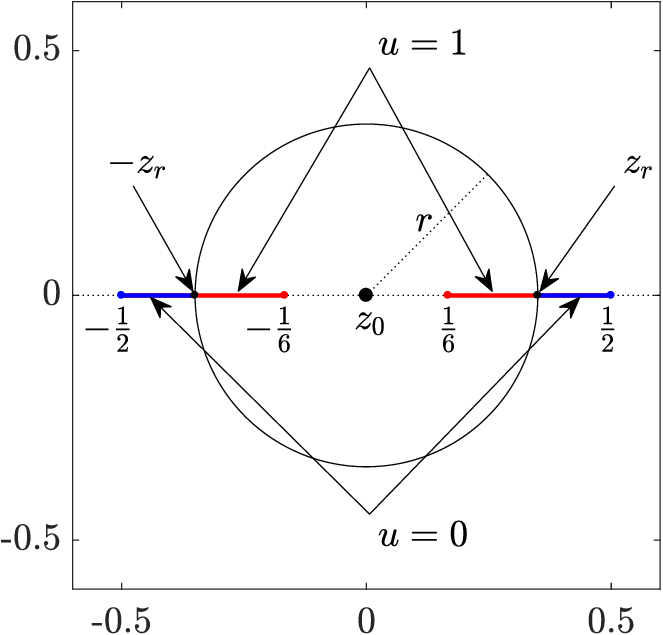}}\hfill
	}
	\caption{The domain $\Omega_2$ and a typical capture circle $\partial B(z_0,r)$ of radius $r$ used in the definition of the $h$-function with basepoint $z_0=-3/2$ (left) and with basepoint $z_0=0$ (right).}
	\label{fig:h30}
\end{figure}

There have been several calculations and analyses of $h$-functions over simply connected domains \cite{SnipWard16}. Recently, the first explicit formulae of $h$-functions in the multiply connected setting were derived by Green \textit{et al.}~\cite{gswc}. Their formulae were constructed by making judicious use of the calculus of the Schottky-Klein prime function, a special transcendental function which is ideally suited to solving problems in multiply connected domains, and the reader is referred to the monograph by Crowdy \cite{CrowdyBook} for further details. 
In~\cite{gswc}, the planar domains of interest were chosen to be a particular class of multiply connected slit domain, i.e. unbounded domains in the complex plane whose boundary components consist of a finite number of finite-length linear segments lying on the real axis. Indeed, this choice of planar domain was motivated in part by the middle-thirds Cantor set. The formulae in \cite{gswc} have been generalized in~\cite{ArMa} for other classes of multiply connected domains and basepoint locations.

Whilst it is important to point out that there is now reliable software available to compute the prime function  based on the analysis of Crowdy \textit{et al.}~\cite{ckgn}, the number of computational operations needed to be executed is relatively high for domains with many boundary components such as those which will be considered in this paper. To this end, to compute $h$-functions for \textit{highly} multiply connected domains, we turn to a boundary integral equation (BIE) method which significantly reduces the number of computational operations required while also maintaining a good level of accuracy. 
The proposed method is based on a BIE with the Neumann kernel (see~\cite{Nas-ETNA,Weg-Nas}). 
The integral equation method that will be employed in this paper has also been used in the new numerical scheme for computing the prime function \cite{ckgn}. The integral equation has found multitudinous application in several areas, including fluid stirrers, vortex dynamics, composite materials, and conformal mappings (see~\cite{Nas-ETNA,NG18,NK} and the references cited therein).
In particular, the integral equation has been used in~\cite{LSN17,Nvm} to compute the logarithmic and conformal capacities for the Cantor set and the Cantor dust.

With the middle-thirds Cantor $\mathcal{C}$ set in mind, the goal herein is to compute approximations to the $h$-functions in domains which effectively imitate the $h$-function associated with $\mathcal{C}$.
As mentioned above, to imitate $\mathcal{C}$ in a practical sense, we again consider multiply connected slit domains as in \cite{gswc}, but significantly increase the number of boundary components by successively removing the middle-third portion of each slit, and then calculate the corresponding $h$-functions.
We also use the computed values of these $h$-functions to study numerically their asymptotic behavior.
The calculation of such $h$-functions, together with their asymptotic features, is motivated by the list of open problems in Snipes \& Ward \cite{SnipWard16}.

Computing the $h$-function requires solving a Dirichlet BVP in the domain $\Omega_m$, which is not an easy task since the boundary of the domain consists of the slits $I_1,\ldots,I_m$ and hence the solution of the BVP admits singularities at the slit endpoints. 
Since the Dirichlet BVP is invariant under conformal mapping, one way to overcome such a difficulty is to map the domain $\Omega_m$ onto another domain with a simpler geometry. 
We therefore appeal to the iterative method presented in~\cite{NG18} to numerically compute an unbounded preimage domain $G_m$ bordered by smooth Jordan curves. 
We will summarize the construction of this preimage domain in \S\ref{sec:pre} below.

Our paper mostly has a potential theoretical flavor and has the following structure. 
In Section 2, we introduce the conformal mapping we use, together with an outline of the BIE with the Neumann kernel. 
When the capture circle passes through the gaps between the slits, the values the $h$-function $h(r)$ takes on are equal to the harmonic measures of the slits enclosed by the capture circle with respect to the domain $\Omega_m$. In this case, with the help of a conformal mapping, the BVP~\eqref{eq:bdv-u} will be transformed into an equivalent problem in an unbounded circular domain $G_m$.
The transformed problem is a particular case of the BVP considered in~\cite{Nvm} and hence it can be solved using the BIE method as in~\cite{Nvm}. This will be considered in Section~\ref{sec:step}. 
We calculate the values of $h$-functions associated with slit domains of connectivity $64$ and $1024$ for the basepoint $z_0=-3/2$; and for connectivity $64$ and $2048$ for the basepoint $z_0=0$. However, when the capture circle intersects one of the slits, the right-hand side of the BVP~\eqref{eq:bdv-u} becomes discontinuous. As in  Section~\ref{sec:step}, we use conformal mapping to transform the slit domain $\Omega_m$ onto a circular domain $G_m$, and hence the BVP~\eqref{eq:bdv-u} will be turned into an equivalent problem in $G_m$ whose right-hand side is also discontinuous. With the help of two M{\"o}bius mappings, the BVP in $G_m$ can be modified so that it will have a continuous right-hand side. 
The BIE method is then used to solve this modified BVP and to compute the values of the $h$-function $h(r)$. This case will be considered in Section~\ref{sec:h-fun}.
In Section 5, we analyze numerically the behavior of the $h$-function $h(r)$ when $r$ is close to $1$ for the basepoint $z_0=-3/2$; and when $r$ is close to $1/6$ for the basepoint $z_0=0$. Section 6 provides our concluding remarks.

\section{The conformal mapping}\label{sec:cm}

Computing the values of the $h$-function requires solving the BVP~\eqref{eq:bdv-u} in the unbounded multiply connected domain $\Omega_m$ bordered by slits.
In this section, we review an iterative conformal mapping method from~\cite{AST13,NG18} which `opens up' the slits $I_k$ and consequently maps the domain $\Omega_m$ onto an unbounded multiply connected circular domain $G_m$. 
The BVP~\eqref{eq:bdv-u} will then be transformed into an equivalent problem in the circular domain $G_m$ where it can be solved in a straightforward manner. Solving the transformed problem in the new circular domain $G_m$ will be discussed in Sections~\ref{sec:step} and~\ref{sec:h-fun} below. 

\subsection{The Neumann kernel}\label{sec:NK}

Let $G_m$ be an unbounded multiply connected circular domain obtained by removing $m$ non-overlapping disks from the extended complex plane $\overline{\CC}$. The boundaries of these disks are the circles $C_j$ with centers $c_j$ and radii $r_j$, $j=1,\ldots,m$. 
We parametrize each circle $C_j$ by 
\begin{equation}\label{eq:eta-j}
	\eta_j(t)=c_j+r_j e^{-\i t}, \quad t\in J_j=[0,2\pi], \quad j=1,\ldots,m.
\end{equation}
We define the total parameter domain $J$ as the disjoint union of the $m$ intervals $J_j=[0,2\pi]$, $j=1,\ldots,m$. Thus, the whole boundary $\Gamma$ is parametrized by
\begin{equation}\label{eq:eta}
	\eta(t)= \left\{ \begin{array}{l@{\hspace{0.5cm}}l}
		\eta_1(t),&t\in J_1, \\
		\quad\vdots & \\
		\eta_m(t),&t\in J_m.
	\end{array}
	\right.
\end{equation}
See~\cite{Nas-ETNA,NG18} for more details.

The Neumann kernel $N(s,t)$ is defined for $(s,t)\in J\times J$ by
\begin{equation}\label{eq:N}
	N(s,t) =
	\frac{1}{\pi}\Im\left(\frac{\eta'(t)}{\eta(t)-\eta(s)}\right).
\end{equation}
It is a particular case of the generalized Neumann kernel considered in~\cite{Weg-Nas} when $A(t)=1$, $t\in J$.
We also define the following kernel 
\begin{equation}\label{eq:M}
	M(s,t) =
	\frac{1}{\pi}\Re\left(\frac{\eta'(t)}{\eta(t)-\eta(s)}\right), \quad (s,t)\in J\times J,
\end{equation}
which is a particular case of the kernel $M$ considered in~\cite{Weg-Nas} when $A(t)=1$.
The kernel $N(s,t)$ is continuous and the kernel $M(s,t)$ is singular.
Hence, the integral operator 
\[
\bN\mu(s) = \int_J N(s,t) \mu(t) dt, \quad s\in J,
\]
is compact and the integral operator 
\[
\bM\mu(s) = \int_J  M(s,t) \mu(t) dt, \quad s\in J,
\]
is singular. Further details can be found in~\cite{Weg-Nas}.

\subsection{The preimage domain}\label{sec:pre}

Let $\Omega_m$ be the multiply connected slit domain obtained by removing the $m$ horizontal slits $I_1,\ldots,I_m$, described in Section~\ref{sec:int}, from the extended complex plane $\overline{\CC}$. In this subsection, we will summarize  an iterative method from~\cite{AST13,NG18} for the construction of a preimage unbounded circular domain $G_m$. The method has been tested numerically in several works~\cite{LSN17,NK,Nvm}.

For a fixed $\ell=0,1,2,\ldots$, the closed set $E_\ell$, defined in~\eqref{eq:Ek}, consists of $m=2^\ell$ slits, $I_1,\ldots,I_m$, of equal length $L=(1/3)^\ell$ and with centers $w_j$ for $j=1,\ldots,m$. Our objective now is to find the centers $c_j$ and radii $r_j$ of non-overlapping circles $C_j$ and a conformal mapping $z=F(\zeta)$ from the circular domain $G_m$ exterior to $C=\cup_{j=1}^m C_j$ onto $\Omega_m$. With the normalization
\[
F(\infty)=\infty, \quad \lim_{\zeta\to\infty}(F(\zeta)-\zeta)=0,
\]
such a conformal mapping is unique.

The conformal mapping $z=F(\zeta)$ can be computed using the following boundary integral equation method from~\cite{Nas-Siam1}. Let the function $\gamma$ be defined by
\begin{equation}\label{eq:gam}
	\gamma(t)=\Im\eta(t), \quad t\in J.
\end{equation}
Let also $\mu$ be the unique solution of the boundary integral equation with the Neumann kernel
\begin{equation}\label{eq:ie-g}
	(\bI-\bN)\mu=-\bM\gamma,
\end{equation}
and let the piecewise constant function $\nu=(\nu_1,\nu_2,\ldots,\nu_\ell)$ be given by
\begin{equation}\label{eq:h-g}
	\nu=\left(\bM\mu-(\bI-\bN)\gamma\right)/2,
\end{equation}
i.e., the function $\nu$ is constant on each boundary component $C_j$ and its value on $C_j$ is a real constant $\nu_j$.
Then the function $f$ with the boundary values
\begin{equation}\label{eq:f-rec}
	f(\eta(t))=\gamma(t)+\nu(t)+\i\mu(t)
\end{equation}
is analytic in $G_m$ with $f(\infty)=0$, and the conformal mapping $F$ is 
given by
\begin{equation}\label{eqn:omega-app}
	F(\zeta)=\zeta-\i f(\zeta), \quad \zeta\in G_m\cup\Gamma.
\end{equation}

The application of this method requires that the domain $G_m$ is known. However, in our case, the slit domain $\Omega_m$ is known and the domain $G_m$ is unknown and needs to be determined alongside the conformal mapping $z=F(\zeta)$ from $G_m$ onto $\Omega_m$. 
This preimage domain $G_m$ as well as the conformal mapping $z=F(\zeta)$ will be computed using the iterative method presented in~\cite{NG18} (see also~\cite{AST13}). 
For the convenience of the reader and for the completeness of this paper, we now briefly review a slightly modified version of this iterative method. 
In~\cite{NG18}, the boundary components of the preimage domain $G_m$ are assumed to be ellipses. Here, the slits $I_{j}$, $j=1,2,\ldots,m$, are well-separated and hence the boundary components of the preimage domain $G_m$ are assumed to be circles. The method generates a sequence of multiply connected circular domains $G_m^0,G_m^1,G_m^2,\ldots,$ which converge numerically to the required preimage domain $G_m$. 
Recall that the length of each slit $I_{j}$ is $L$ and the center of $I_{j}$ is $w_j$ for $j=1,\ldots,m$. 
In the iteration step $i=0,1,2,\ldots$, we assume that $G_m^{(i)}$ is bordered by the $m$ circles $C^{(i)}_1,\ldots,C^{(i)}_m$ parametrized by
\begin{equation}\label{eq:eta-i}
	\eta^{(i)}_j(t)=c^{(i)}_j+r^{(i)}_je^{\i t}, \quad 0\le t\le 2\pi, \quad j=1,\ldots,m.
\end{equation}
The centers $c^{(i)}_j$ and radii $r^{(i)}_j$ are computed using the following iterative method:
\begin{enumerate}
	\item Set
	\[
	c^{(0)}_j=w_j, \quad r^{(0)}_j=\frac{\,L\,}{2}, \quad j=1,\ldots,m.
	\]
	\item For $i=1,2,3,\ldots,$
	\begin{itemize}
		\item Compute the conformal mapping from the preimage domain $G^{(i-1)}$ to a canonical horizontal rectilinear slit domain $\Omega^{(i)}$ which is the entire $\zeta$-plane with $m$ horizontal slits $I^{(i)}_{j}$, $j=1,\ldots,m$ (using the method presented in equations~\eqref{eq:gam}--\eqref{eqn:omega-app} above). Let $L^{(i)}_j$ denote the length of the slit $I^{(i)}_{j}$ and let $w^{(i)}_j$ denote its center.
		\item Define	
		\[
		c^{(i)}_j = c^{(i-1)}_j-(w^{(i)}_j-w_j), \quad 
		r^{(i)}_j = r^{(i-1)}_j-\frac{\,1\,}{4}(L^{(i)}_j -L), \quad j=1,\ldots,m.
		\]

	\end{itemize}
	\item Stop the iteration if 
	\[
	\frac{1}{2m}\sum_{j=1}^{m}\left(|w^{(i)}_j - w_j|+|L^{(i)}_j -L|\right)<\varepsilon \quad{\rm or}\quad i>{\tt Max}
	\]
	where $\varepsilon$ is a given tolerance and ${\tt Max}$ is the maximum number of iterations allowed.		
\end{enumerate}

The above iterative method generates sequences of parameters $c^{(i)}_j$ and $r^{(i)}_j$ that converge numerically to $c_j$ and $r_j$, respectively, and then the boundary components of the preimage domain $G_m$ are parametrized by~\eqref{eq:eta-j}. In our numerical implementations, we used $\varepsilon=10^{-14}$ and ${\tt Max}=100$.

It is clear that in each iteration of the above method, it is required to solve the integral equation with the Neumann kernel~(\ref{eq:ie-g}) and to compute the function $\nu$ in~(\ref{eq:h-g}).
In this paper, we approximate the solution of the integral equation~(\ref{eq:ie-g}) and the function $\nu$ using the fast numerical method presented in~\cite{Nas-ETNA}. We shall now briefly describe the implementation of this numerical method and refer the reader to~\cite{Nas-ETNA} for the details (see also~\cite{LSN17,NG18}). 

Since the boundary components $C_1,\ldots,C_m$ of the domain $G_m$ are circles and hence $C^\infty$ smooth boundaries, the boundary integral equation~\eqref{eq:ie-g} can be solved accurately by the Nystr\"om method with the trapezoidal rule~\cite{Atk97}. This leads to an $mn\times mn$ linear system with a dense non-symmetric coefficient matrix, where $n$ is the number of nodes in the discretization 	of each boundary component. 
This system can be solved using the GMRES iterative method. Each step of this method requires one multiplication with the coefficient matrix which can be computed efficiently in $\mathcal{O}(m n)$ operations using the Fast Multipole Method (FMM).
The number of GMRES iterations for obtaining a very good approximation to the exact solution is virtually independent of the given domain and the number of nodes in the discretization of its boundary.
This method for solving the integral equation~\eqref{eq:ie-g} for $\mu$ and then calculating $\nu$ in~\eqref{eq:h-g} has been implemented
in the {\sc Matlab} function {\tt fbie} shown in~\cite[Figure~4.1]{Nas-ETNA}.
The function {\tt fbie} uses {\sc Matlab}'s built-in {\tt gmres} function together with the function {\tt zfmm2dpart} from the fast multipole toolbox
FMMLIB2D~\cite{Gre-Gim12}. The main inputs of the function {\tt fbie} are the discretizations of the function $\eta(t)$ given by~\eqref{eq:eta}, its derivative $\eta'(t)$, and the function $\gamma(t)$ defined in terms of $\eta(t)$ by~\eqref{eq:gam}.

To implement {\tt fbie}, we discretize the interval
$[0, 2\pi]$ by the $n$ equidistant nodes
\begin{equation}\label{eq:s_i}
	s_k = (k-1) \frac{2 \pi}{n}, \quad k = 1, \ldots, n,
\end{equation}
and write $\bs = [s_1, \ldots, s_n]$ where $n$ is an even integer.  
Then $\bt = [\bs, \bs, \ldots, \bs]^T$, which consists of $m$ copies of $\bs$, provides a discretization of the total parameter domain $J$. This leads to the following discretizations
\begin{equation}\label{eq:dis_s}
	\eta(\bt) = [\eta_1(\bs), \eta_2(\bs), \ldots, \eta_\ell(\bs)]^T, \quad
	\eta'(\bt), \quad
	\gamma(\bt) =\Im\eta(\bt),
\end{equation}
of the functions $\eta(t)$, $\eta'(t)$, and $\gamma(t)$. 
We store these discretized functions in the $mn\times1$ vectors \texttt{et}, \texttt{etp}, \texttt{gam}, respectively, and call
\begin{verbatim}
	[mu,nu] = FBIE(et,etp,ones(size(et)),gam,n,5,[],tol,maxit).
\end{verbatim}
Here, {\tt []} means that GMRES is used without restart, {\tt tol} is the convergence criterion used within GMRES, and {\tt maxit}
is the maximal number of GMRES iterations. In the numerical experiments presented in this paper, we used {\tt tol=1e-13} and {\tt maxit=100}.
The boundary components of $\Omega_m$, i.e., the $m$ slits $I_1,\ldots,I_m$, are well-separated and hence the circles $C_1,\ldots,C_m$, which are the boundary components of the domain $G_m$, are also well-separated (see Figure~\ref{fig:map} for $m=4$). Thus, accurate results can be obtained even with small values of $n$~\cite{Nvm}. In our numerical computations, we used $n=16$. For the {\sc Matlab} function \verb|zfmm2dpart|, we chose the FMM precision flag \verb|iprec=5| which means that the tolerance of the FMM method is $0.5\times10^{-15}$. The complexity of solving this integral equation, and hence the complexity of each iteration in the above iterative method, is $\mathcal{O}(mn\log n)$ operations.

The output of {\tt fbie} are the $mn\times1$ vectors {\tt mu} and {\tt nu} which are the approximations of the discretizations $\mu(\bt)$ and $\nu(\bt)$ of the functions $\mu(t)$ and $\nu(t)$, respectively. It then follows from~\eqref{eq:f-rec} that  approximate discretizations of the boundary values of the analytic function $f(\zeta)$ are given by
\[
f(\eta(\bt))=\gamma(\bt)+\nu(\bt)+\i\mu(\bt).
\]
The values of $f(\zeta)$ for $\zeta\in G_m$ can then be computed using the Cauchy integral formula. Finally, the values of the conformal mapping $z=F(\zeta)$ are computed for $\zeta\in G_m\cup\Gamma$ via~\eqref{eqn:omega-app}.

\section{Step heights of the $h$-function}\label{sec:step}

In this section, we determine the step heights of the $h$-functions of interest. We refer to the `step heights' as the constant values of $h(r)$ coinciding with when the capture circle passes through the gaps between the slits, for the two basepoint locations $z_0=-3/2$ and $z_0=0$ (see Figure~\ref{fig:h30}). The restriction of the $h$-function $h(r)$ to these values of $r$ will be denoted by $\omega(r)$. The step heights of the $h$-function for an example with two slits have been computed in~\cite{SnipWard16} and for two, four and eight slits in~\cite{gswc}. Furthermore, it should be pointed out that there are several methods for computing harmonic measures in multiply connected domains, each of which may alternatively be used to compute the step heights of the $h$-function; see e.g.,~\cite{CrowdyBook,Del-harm,Tre-Gre,garmar,Tre-Ser} and the references cited therein.

\begin{figure}[ht] %
	\centerline{\hfill
		\scalebox{0.4}{\includegraphics[trim=0 0 0 0,clip]{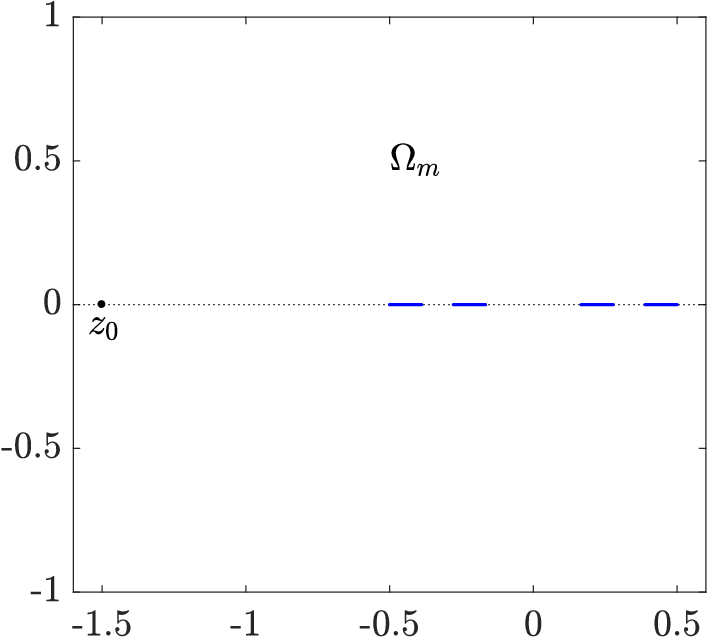}}\hfill
		\scalebox{0.4}{\includegraphics[trim=0 0 0 0,clip]{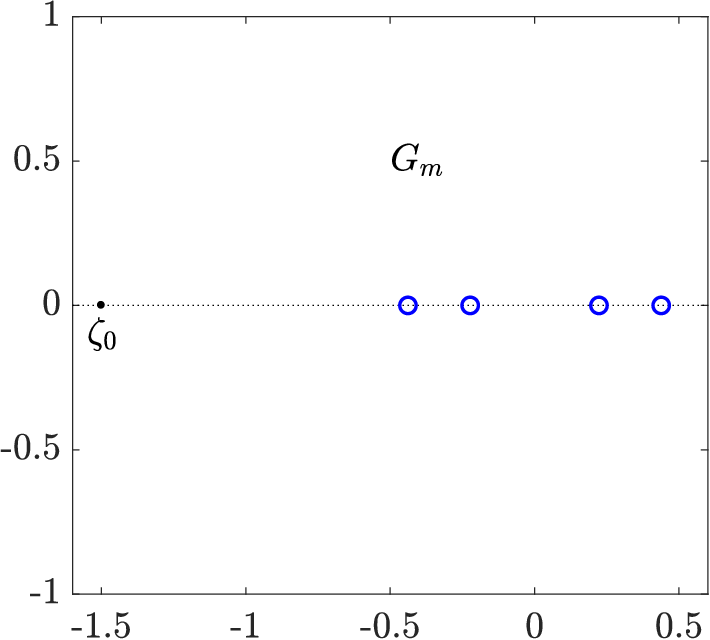}}\hfill
	}
	\caption{The slit domain $\Omega_m$ (left) and the circular domain $G_m$ (right) with basepoint $z_0=-3/2$ in the case when $m=4$.}
	\label{fig:map}
\end{figure}

By the mapping function $\zeta=F^{-1}(z)$ described in \S\ref{sec:pre}, the multiply connected slit domain $\Omega_m$ (exterior to the $m$ slits $I_j$, $j=1,\ldots,m$) is mapped conformally onto a multiply connected circular domain $G_m$ exterior to $m$ disks such that each slit $I_{j}$ is mapped onto the circle $C_j$ with the center $c_j$ and radius $r_j$, $j=1,\ldots,m$. The basepoint $z_0$ is also mapped by $\zeta=F^{-1}(z)$ onto a point $\zeta_0$ in the domain $G_m$.
Owing to the symmetry of $\Omega_m$ with respect to the real axis, the centers of the disks as well as the point $\zeta_0$ are on the real line (see Figure~\ref{fig:map} for $m=4$ and $z_0=-3/2$).

\subsection{Harmonic measures}

For $k=1,\ldots,m$, let $\sigma_k$ be the harmonic measure of $C_k$ with respect to $G_m$, i.e., $\sigma_k(\zeta)$ is the unique solution of the Dirichlet problem:
\begin{subequations}\label{eq:bdv-sig}
	\begin{align}
		\label{eq:sig-Lap}
		\nabla^2 \sigma_k(\zeta) &= 0 \quad\quad \mbox{if }\zeta\in G_m, \\
		\label{eq:sig-j}
		\sigma_k(\zeta)&= \delta_{k,j} \quad \mbox{if }\zeta\in C_j, \quad j=1,\ldots,m, 
	\end{align}
\end{subequations}
where $\delta_{k,j}$ is the Kronecker delta function.
The function $\sigma_k$ is assumed to be bounded at infinity. 
The harmonic function $\sigma_k$ is the real part of an analytic function $g_k$ in $G_m$ which is not necessarily single-valued. The function $g_k$ can be written as~\cite{Gak,garmar}
\begin{equation}\label{eq:F-u}
	g_k(\zeta)=b_k+f_k(\zeta)-\sum_{j=1}^{m} a_{kj}\log(\zeta-c_j)
\end{equation}
where  $c_j$ is the center of the circle $C_j$, $f_k$ is a single-valued analytic function in $G_m$ with $f_k(\infty)=0$, $b_k$ and $a_{k,j}$ are undetermined real constants such that $\sum_{j=1}^{m}a_{kj}=0$ for $k=1,\ldots,m$. The condition $\sum_{j=1}^{m}a_{kj}=0$ implies that $g_k(\infty)=b_k$ since $f_k(\infty)=0$. Since we are interested in computing the real function $\sigma_k=\Re[g_k]$, we may assume that $b_k$ is real, and then $b_k=g_k(\infty)=\sigma_k(\infty)$.

The BVP~\eqref{eq:bdv-sig} above is a particular case of the problem considered in~\cite[Eq.~(4)]{Nvm}; hence, it can be solved by the method presented in~\cite{Nvm} which is reviewed in the following subsection.

\subsection{Computing the harmonic measures}

For $k=1,\ldots,m$, it follows from \eqref{eq:F-u} that computing the harmonic measure $\sigma_k=\Re[g_k]$ requires computing the values of the analytic function $f_k$ and the values of the $m+1$ real constants $a_{k,1},\ldots,a_{k,m},b_k$. The constants $a_{k,1},\ldots,a_{k,m}$ in~\eqref{eq:F-u} can be computed as described in Theorem~3 in~\cite{Nvm} (with $\ell=0$ and $A=1$). 
For each $j=1,2,\ldots,m$, let the function $\gamma_j$ be defined by
\begin{equation}\label{eq:gam-j}
	\gamma_j(t)=\log|\eta(t)-c_j|,
\end{equation}
let $\mu_j$ be the unique solution of the boundary integral equation with the Neumann kernel
\begin{equation}\label{eq:ie}
	(\bI-\bN)\mu_j=-\bM\gamma_j,
\end{equation}
and let the piecewise constant function $\nu_j=(\nu_{1,j},\nu_{2,j},\ldots,\nu_{m,j})$ be given by
\begin{equation}\label{eq:hj}
	\nu_j=\left(\bM\mu_j-(\bI-\bN)\gamma_j\right)/2.
\end{equation}
Then, for each $k=1,2,\ldots,m$, the boundary values of the function $f_k(\zeta)$ in~\eqref{eq:F-u} are given by
\begin{equation}\label{eq:fk}
	f_k(\eta(t))=\sum_{j=1}^m a_{kj}\left(\gamma_j(t)+\nu_j(t)+\i\mu_j(t)\right)
\end{equation}
and the $m+1$ unknown real constants $a_{k,1},\ldots,a_{k,m},b_k$ are the unique solution of the linear system
\begin{equation}\label{eq:sys-method}
	\left[\begin{array}{ccccc}
		\nu_{1,1}    &\nu_{1,2}    &\cdots &\nu_{1,m}      &1       \\
		\nu_{2,1}    &\nu_{2,2}    &\cdots &\nu_{2,m}      &1       \\
		\vdots     &\vdots     &\ddots &\vdots       &\vdots  \\
		\nu_{m,1}    &\nu_{m,2}    &\cdots &\nu_{m,m}      &1       \\
		1          &1          &\cdots &1            &0       \\
	\end{array}\right]
	\left[\begin{array}{c}
		a_{k,1}    \\a_{k,2}    \\ \vdots \\ a_{k,m} \\  b_k 
	\end{array}\right]
	= \left[\begin{array}{c}
		\delta_{k,1} \\  \delta_{k,2} \\  \vdots \\ \delta_{k,m} \\ 0  
	\end{array}\right].
\end{equation}
The integral operators $\bN$ and $\bM$ in~\eqref{eq:ie} and~\eqref{eq:hj} are the same operators introduced in Section~\ref{sec:NK}.

It is clear that computing the $m+1$ unknown real constants $a_{k,1},\ldots,a_{k,m},b_k$ requires solving $m$ integral equations with the Neumann kernel~\eqref{eq:ie} and computing $m$ piecewise constant functions $\nu_j$ in~\eqref{eq:hj} for $j=1,\ldots,m$. 
This can be done using the {\sc Matlab} function \verb|fbie| as described in Section~\ref{sec:pre}. 
The complexity of solving each of these integral equations is $\mathcal{O}(mn\log n)$ operations and hence solving the $m$ integral equations~\eqref{eq:ie} requires $\mathcal{O}(m^2n\log n)$ operations.

The linear system~\eqref{eq:sys-method} has an $(m+1)\times(m+1)$ constant coefficient matrix and $m$ different right-hand sides. 
These $m$ linear systems are solved in $\mathcal{O}(m^3)$ operations by computing the inverse of the coefficient matrix and then using this inverse matrix to compute the solution for each right-hand side. In doing so, we obtain the values of the real constants $a_{k,1},\ldots,a_{k,m},b_k$. The boundary values of the analytic function $f_k$ are given by~\eqref{eq:fk}, and hence its values in the domain $G_m$ can be computed by the Cauchy integral formula. The harmonic measure $\sigma_k$ can then be computed for $\zeta\in G_m$ by
\begin{equation}\label{eq:sgj-gj}
	\sigma_k(\zeta)=\Re g_k(\zeta)
\end{equation}
where $g_k(\zeta)$ is given by~\eqref{eq:F-u}.

\subsection{Basepoint $z_0=-3/2$}

The capture circle of radius $r$ intersects with the real line at two real numbers where the largest of these numbers is denoted by $z_r$ (see Figure~\ref{fig:h30} (left)).
By `step heights', we mean the constant values of $h(r)$ when $r$ is such that $z_r$ lies in between a pair of slits $I_{j}$ and $I_{j+1}$, $j=1,\ldots,m-1$. In addition, $h(r) = 0$ when $r$ is such that $z_r$ lies strictly to the left of all $m$ slits (i.e., $0\le r\le1$), and $h(r) = 1$ when $r$ is such that $z_r$ lies strictly to the right of all $m$ slits (i.e., $2\le r<\infty$). 

Here, we use the harmonic measures $\sigma_1,\ldots,\sigma_m$ to compute the $m-1$ step heights of the $h$-functions for the slit domains $\Omega_m$ (the exterior domain of the closed set $E_\ell$) with basepoint location $z_0=-3/2$. That is, we compute the values $\omega(r)$ of these $h$-functions at those values of $r$ corresponding to capture circles that pass through the $m-1$ gaps between the slits $I_1,\ldots,I_m$. 
In the next section, we compute the values of the $h$-functions for $\Omega_m$ for the remaining values of $r$, i.e. those values of $r$ which correspond to capture circles intersecting a slit $I_{j}$, for some $j=1,\ldots,m$. 

\begin{figure}[ht] %
	\centerline{
		\scalebox{0.6}{\includegraphics[trim=0 0 0 0,clip]{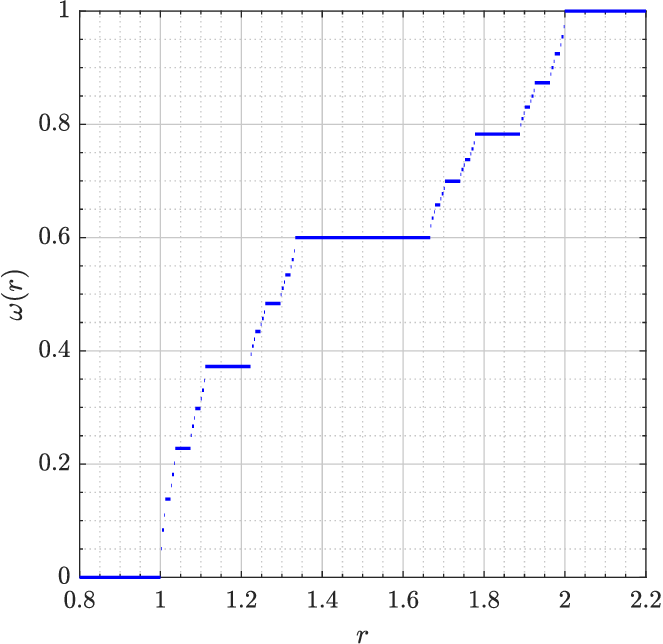}}
		\hfill
		\scalebox{0.6}{\includegraphics[trim=0 0 0 0,clip]{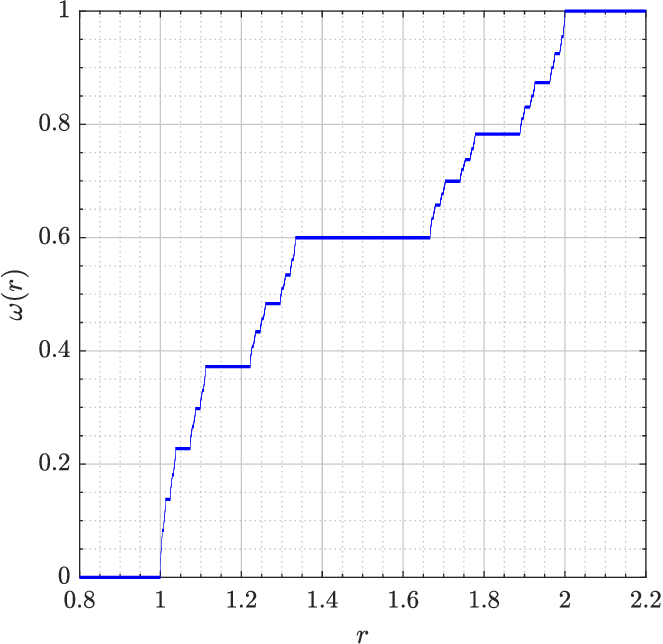}}
	}
	\caption{The step height of the $h$-function for the domain $\Omega_{64}$ (left) and $\Omega_{1024}$ (right) when $z_0=-3/2$.}
	\label{fig:hm-3}
\end{figure}

More precisely, the height $\omega(r)$ of the $k$th step of the $h$-function associated with $\Omega_m$ can be written in terms of the the harmonic measures $\sigma_1,\ldots,\sigma_m$ as 
\begin{equation}\label{eq:wr-3}
	\omega(r)=\sum_{j=1}^{k}\sigma_j(\zeta_0), \quad k=1,2,\ldots,m-1,
\end{equation}
where $\zeta_0=F^{-1}(z_0)$, $\sigma_j$ is the harmonic measure of $C_j$ with respect to $G_m$, and $k$ is the largest integer such that the slit $I_{k}$ is inside the circle $|z-z_0|=r$ and the slit $I_{k+1}$ is outside the circle $|z-z_0|=r$.
The computed numerical values for the step heights of the $h$-function for the domains $\Omega_2$, $\Omega_4$, and $\Omega_8$ are presented in Table~\ref{tab:3}. These values agree with the values computed by the method presented in~\cite{gswc}. The plots of the step height as a function of $r$ for the domains $\Omega_{64}$ and $\Omega_{1024}$ are shown in Figure~\ref{fig:hm-3}.

\begin{table}[h]
	\caption{The $m-1$ step heights of the $h$-function for the domains $\Omega_m$ when $z_0=-3/2$ for $m=2,4,8$, computed using the proposed method and the method presented in~\cite{gswc}.}
	\label{tab:3}
	\centering
	\begin{tabular}{lll|lll}  \hline
		\multicolumn{3}{c|}{The proposed method}  & \multicolumn{3}{c}{Green \emph{et al.}~\cite{gswc}}  \\ \hline
		&              & $0.23081722$  &              &              & $0.23081722$  \\ 
		& $0.37725094$ & $0.37469279$  &              & $0.37725094$ & $0.37469279$ \\
		&              & $0.48515843$  &              &              & $0.48515843$ \\
		$0.60527819$ & $0.60254652$ & $0.60117033$  & $0.60527819$ & $0.60254652$ & $0.60117033$ \\
		&              & $0.70056784$  &              &              & $0.70056784$ \\
		& $0.78306819$ & $0.78288753$  &              & $0.78306819$ & $0.78288753$ \\
		&              & $0.87276904$  &              &              & $0.87276904$ \\
		\hline
	\end{tabular}
\end{table}

\subsection{Basepoint $z_0=0$}

We now compute the step heights of the $h$-functions for $\Omega_m$ with basepoint location $z_0=0$. That is, we compute the values of these $h$-functions at those values of $r$ corresponding to capture circles that pass through gaps between the slits $I_{k}$ for $k=1,\ldots,m$.	
In this case, the capture circle of radius $r$ intersects with the real line at two real numbers $\pm z_r$ where $z_r>0$ (see Figure~\ref{fig:h30} (right)).
Here, $h(r) = 0$ when $r$ is such that $z_r$ lies strictly in-between the two slits $I_{m/2}$ and $I_{m/2+1}$ (i.e., $0\le r\le1/6$), and $h(r) = 1$ when $r$ is such that $z_r$ lies strictly to the right of all $m$ slits (i.e., $1/2\le  r<\infty$). 

\begin{figure}[ht] %
	\centerline{
		\scalebox{0.6}{\includegraphics[trim=0 0 0 0,clip]{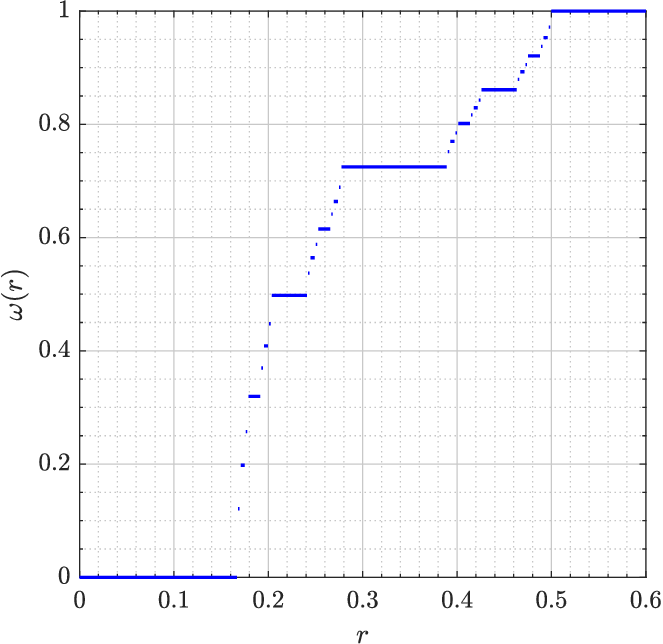}}
		\hfill
		\scalebox{0.6}{\includegraphics[trim=0 0 0 0,clip]{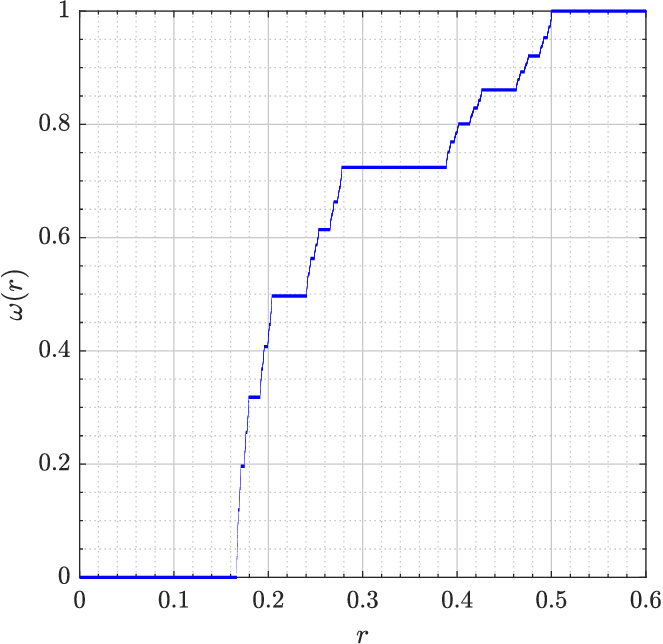}}}
	\caption{The step height of the $h$-function for the domain $\Omega_{64}$ (left) and $\Omega_{2048}$ (right) when $z_0=0$.}
	\label{fig:hm-0}
\end{figure}

For this case, the height $\omega(r)$ of the $k$th step of the $h$-function associated with $\Omega_m$ is thus
\begin{equation}\label{eq:wr-0}
	\omega(r)=\sum_{j=1}^{k}\left(\sigma_{\frac{m}{2}-j+1}(\zeta_0)+\sigma_{\frac{m}{2}+j}(\zeta_0)\right)
\end{equation}
where $\zeta_0=F^{-1}(z_0)$, $\sigma_j$ is the harmonic measure of $C_j$ with respect to $G_m$, and $k$ is the largest integer such that the two slits $I_{m/2-j+1}$ and $I_{m/2+j}$ are inside the circle $|z-z_0|=r$ and the two slits $I_{m/2-j}$ and $I_{m/2+j+1}$ are outside the circle $|z-z_0|=r$.
The computed numerical values for the step heights of the $h$-function for the domains $\Omega_4$, $\Omega_8$, and $\Omega_{16}$ are presented in Table~\ref{tab:0}. The plots of the step height as a function of $r$ for the domains $\Omega_{64}$ and $\Omega_{2048}$ are shown in Figure~\ref{fig:hm-0}.

\begin{table}[h]
	\caption{The $m/2-1$ step heights of the $h$-function for the domains $\Omega_m$ when $z_0=0$ for $m=4,8,16$.}
	\label{tab:0}
	\centering
	\begin{tabular}{ccc}
		\hline
		$\Omega_4$       & $\Omega_8$        & $\Omega_{16}$  \\ \hline
		&                   & $0.32412730$  \\
		& $0.50657767$      & $0.50171513$  \\
		&                   & $0.61794802$  \\
		$0.73555154$      & $0.72992958$      & $0.72702487$  \\
		&                   & $0.80333761$  \\
		& $0.86334249$      & $0.86206402$  \\
		&                   & $0.92098300$  \\
		\hline
	\end{tabular}
\end{table}

\section{Numerical computation of the $h$-function}\label{sec:h-fun}

The step heights $\omega(r)$ of the $h$-function, i.e., the values of the $h$-function $h(r)$ coinciding with when the capture circle passes through the gaps between the slits $I_k$, $k=1,\ldots,m$, have been computed in the previous section. In this section, we present a method for computing the values of the $h$-function $h(r)$ when the capture circle intercepts with the slits $I_k$ for the two cases of basepoints considered in the previous section.

\subsection{Basepoint $z_0=-3/2$}

\begin{figure}[ht] %
	\centerline{
		\scalebox{0.4}{\includegraphics[trim=0 0 0 0,clip]{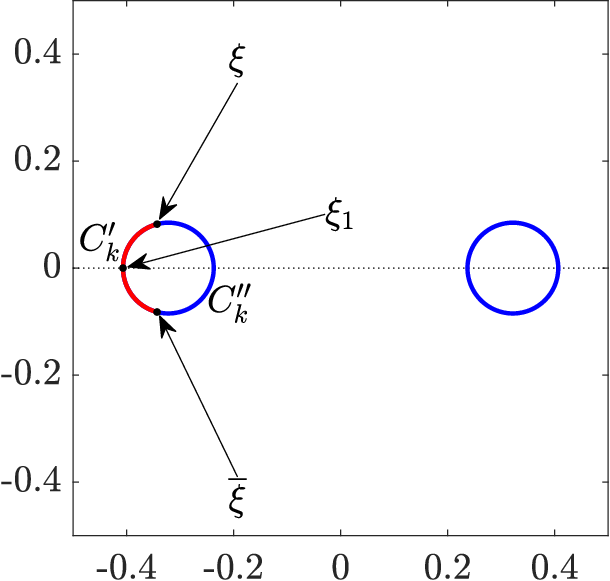}}
	}
	\caption{The arcs $C'_k$ and $C''_k$ for $z_0=-3/2$.}
	\label{fig:h3d}
\end{figure}

Assume that the capture circle intercepts the slit $I_k$ for some $k=1,...,m$. 
Owing to the up-down symmetry of the domain $\Omega_m$ in the $z$-plane, we can assume that the pair of preimages of the intersection point of the capture circle with the slit $I_k$ are the points $\xi$ and $\overline{\xi}$ on the circle $C_k$.
Let us also assume that $\xi_1$ is the intersection of the circle $C_k$ with the real axis lying on the arc between $\overline{\xi}$ and $\xi$ that is closest to the point $\zeta_0$. This arc will be denoted by $C'_k$. The remaining part of the circle is denoted by $C''_k$ (see Figure~\ref{fig:h3d}). 
Let the function $U_k(\zeta)$ be the unique solution of the Dirichlet problem:
\begin{subequations}\label{eq:bdv-U3}
	\begin{align}
		\label{eq:U3-Lap}
		\nabla^2 U_k(\zeta) &= 0 \quad \mbox{if }\zeta\in G_m, \\
		\label{eq:U3-m}
		U_k(\zeta)&= 1 \quad \mbox{if }\zeta\in C_j, \quad j=1,\ldots,k-1, \\
		\label{eq:U3-k'}
		U_k(\zeta)&= 1 \quad \mbox{if }\zeta\in C'_k,  \\
		\label{eq:U3-k''}
		U_k(\zeta)&= 0 \quad \mbox{if }\zeta\in C''_k,  \\
		\label{eq:U3-p}
		U_k(\zeta)&= 0 \quad \mbox{if }\zeta\in C_j, \quad j=k+1,\ldots,m, 
	\end{align}
\end{subequations}
where the function $U_k$ is assumed to be bounded at infinity. 

Note that the boundary conditions on the circle $C_k$ are not continuous. To remove this discontinuity, let us now introduce a function $\Psi_k$ in terms of the two M\"obius maps
\begin{equation}
	\psi(\zeta,\xi,\xi_1)=\frac{(\zeta-\xi)(\xi_1-\overline{\xi})+\i(\zeta-\overline{\xi})(\xi_1-\xi)}{(\zeta-\overline{\xi})(\xi_1-\xi)+\i(\zeta-\xi)(\xi_1-\overline{\xi})}, \quad
	\phi(w)=\frac{w-\i}{\i w-1}.
\end{equation}
Note that $\psi(\zeta,\xi,\xi_1)$ maps the exterior of the unit disc onto the unit disk such that the three points $\overline{\xi}$, $\xi_1$, and $\xi$ on the unit circle are mapped onto the three points $-\i$, $1$, and $\i$, respectively, and $\phi(w)$ maps the unit disc onto the upper half-plane such that the three points $-\i$, $1$, and $\i$ are mapped onto the three points $\infty$, $-1$, and $0$, respectively.
Then, for $\xi\in C_k$, we define 
\begin{equation}\label{eq:Psi}
	\Psi_k(\zeta)=\frac{1}{\pi}\Im \log \phi(\psi(\zeta,\xi,\xi_1)), \quad \zeta\in G_m.
\end{equation}
This function $\Psi_k(\zeta)$ is harmonic everywhere in the domain $G_m$ exterior to the $m$ discs, bounded at infinity, equal to 1 on the arc $C'_k$, and equal to $0$ on the arc $C''_k$.
Hence, for $k=1,\ldots,m$, the function 
\begin{equation}\label{eq:Uk-Vk}
	V_k(\zeta)=U_k(\zeta)-\Psi_k(\zeta)
\end{equation}
is bounded at infinity and satisfies the following Dirichlet problem:
\begin{subequations}\label{eq:bdv-V3}
	\begin{align}
		\label{eq:V3-Lap}
		\nabla^2 V_k(\zeta) &= 0  ~~~~~~~~~~~~~~\,    \mbox{if }\zeta\in G_m, \\
		\label{eq:V3-m}
		V_k(\zeta)&= 1-\Psi_k(\zeta) \quad   \mbox{if }\zeta\in C_j, \quad j=1,\ldots,k-1, \\
		\label{eq:V3-k'}
		V_k(\zeta)&= 0         ~~~~~~~~~~~~~~\,        \mbox{if }\zeta\in C_k,  \\
		\label{eq:V3-p}
		V_k(\zeta)&= -\Psi_k(\zeta)  ~~~~~~  \mbox{if }\zeta\in C_j, \quad j=k+1,\ldots,m. 
	\end{align}
\end{subequations}
The boundary conditions of the new problem~\eqref{eq:bdv-V3} are now continuous on all circles.

\begin{figure}[htb] %
	\centerline{\hfill
		\scalebox{0.4}{\includegraphics[trim=0 0 0 0,clip]{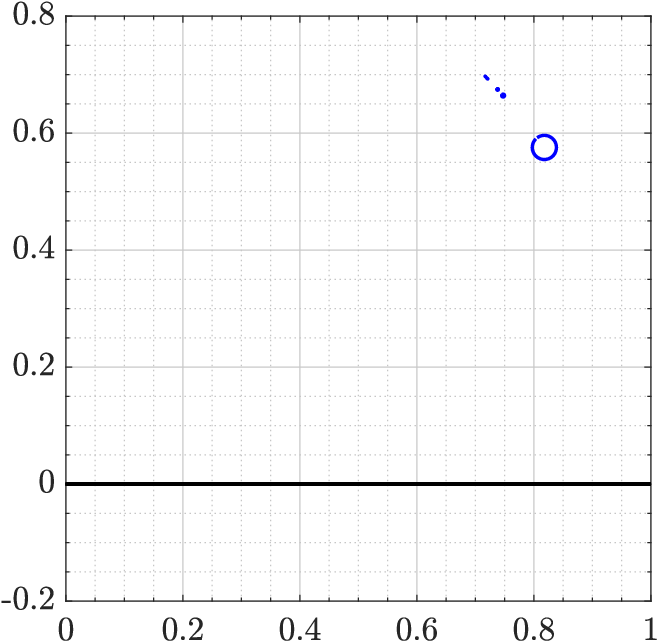}}
		\hfill
		\scalebox{0.4}{\includegraphics[trim=0 0 0 0,clip]{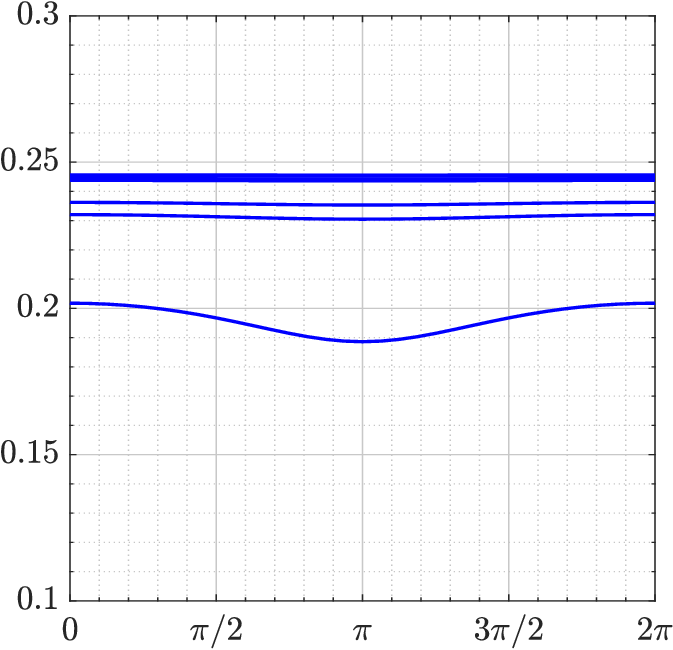}}\hfill
	}
	\caption{The image of the circles $C_j$ (left) and the values of the functions in the right-hand side of the boundary conditions on~\eqref{eq:bdv-V3} (right) for $m=8$, $k=1$, and $\xi=c_k+r_k e^{3\pi\i/4}\in C_k$.}
	\label{fig:small}
\end{figure}

Note that $\phi(\psi(\zeta,\xi,\xi_1))$ is itself a M\"obius map which maps the circle $C_k$ onto the real line. For $j=1,\ldots,m$, $j\ne k$, it  maps the circle $C_j$ onto a circle of very small radius in the upper half-plane (see Figure~\ref{fig:small} (left)). Hence in problem~\eqref{eq:bdv-V3}, the values of the function on the right-hand side of the boundary conditions are almost constant (see Figure~\ref{fig:small} (right)). As such, and in view of (\ref{eq:Psi}), we may approximate
\begin{equation}
	\Psi_k(\zeta) \approx P_{kj}= \Psi_k(c_j), \quad \zeta \in C_j, \quad j=1,\ldots,m,
\end{equation} 
where $c_j$ is the center of the circle $C_j$. Thus, the Dirichlet problem~\eqref{eq:bdv-V3} becomes
\begin{subequations}\label{eq:bdv-V23}
	\begin{align}
		\label{eq:V23-Lap}
		\nabla^2 V_k(\zeta) &= 0 ~~~~~~~~~~\:\: \mbox{if }\zeta\in G_m, \\
		\label{eq:V23-m}
		V_k(\zeta)&= 1-P_{kj} \quad \mbox{if }\zeta\in C_j, \quad j=1,\ldots,k-1, \\
		\label{eq:V23-k'}
		V_k(\zeta)&= 0 ~~~~~~~~~~\:\: \mbox{if }\zeta\in C_k,  \\
		\label{eq:V23-p}
		V_k(\zeta)&= -P_{kj}  ~~~~\:\:\: \mbox{if }\zeta\in C_j, \quad j=k+1,\ldots,m. 
	\end{align}
\end{subequations}
The function $V_k$ can be then written in terms of the harmonic measures $\sigma_1,\ldots,\sigma_m$ through
\[
V_k(\zeta)=\sum_{j=1}^{k-1} (1-P_{kj})\sigma_j(\zeta)-\sum_{j=k+1}^{m} P_{kj}\sigma_j(\zeta).
\]
Thus, by~\eqref{eq:Uk-Vk}, the unique solution $U_k(\zeta)$ to the Dirichlet problem~\eqref{eq:bdv-U3} is given by
\begin{equation}\label{eq:U_k}
	U_k(\zeta) = \Psi_k(\zeta)+\sum_{j=1}^{k-1} \sigma_j(\zeta) -\sum_{\begin{subarray}{c} j=1\\j\ne k\end{subarray}}^{m} P_{kj}\sigma_j(\zeta).
\end{equation}

\begin{figure}[htb] %
	\centerline{\hfill
		\scalebox{0.4}{\includegraphics[trim=0 0 0 0,clip]{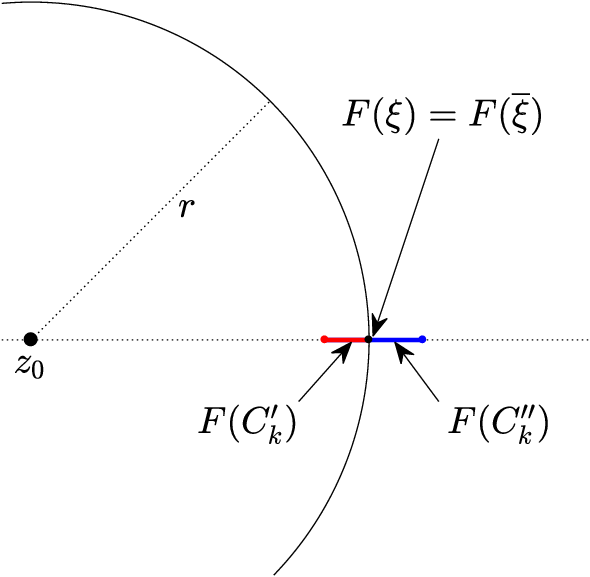}}
		\hfill
		\scalebox{0.4}{\includegraphics[trim=0 0 0 0,clip]{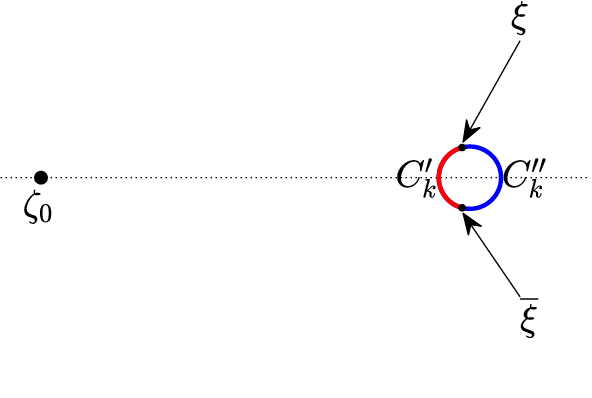}}\hfill
	}
	\caption{Computing the value of $r$ in~\eqref{eq:r-xi}.}
	\label{fig:3hr}
\end{figure}

Due to the symmetry of both domains $\Omega_m$ and $G_m$ with respect to the real line and since the image of the circle $C_k$ under the conformal mapping $F$ is the slit $I_k$, we have $F(\xi)=F(\overline{\xi})\in I_k$ and the image of the arc $C_k'$ is the part of $I_k$ lying to the left of $F(\xi)$ (see Figure~\ref{fig:3hr}). 
For a given $r$ such that the capture circle passes through the slit $I_k$, for $k=1,...,m$, then $\xi$ is on the circle $C_k$ and depends on $r$ as (see Figure~\ref{fig:3hr}) 
\begin{equation}\label{eq:r-xi}
	r= F(\xi)-z_0.
\end{equation}
For this value of $r$, the value of the $h$-function is then given by
\begin{equation}\label{eq:hr3-xi}
	h(r)=U_k(\zeta_0).
\end{equation}

In~\eqref{eq:r-xi}, we may assume that $r$ is given and then solve the non-linear equation~\eqref{eq:r-xi} for $\xi\in C_k$. Alternatively, we may assume that $\xi$ on the circle $C_k$ is given and then compute the values of $r$ through~\eqref{eq:r-xi}, which is simpler. 
For computing the value of $r$ in~\eqref{eq:r-xi}, the value of the mapping function $F(\xi)$ can be approximated numerically for $\xi\in C_k$ as follows.
Since $\xi\in C_k$, then it follows from the parametrization~\eqref{eq:eta-j} of the circle $C_k$ that $\xi=\eta_k(s)$ for some $s\in J_k=[0,2\pi]$. Then, it follows from the method described in Section~\ref{sec:pre} that $F(\xi)=\xi-\i f(\xi)$ where
\[
f(\xi)=f(\eta_k(s))=\gamma_k(s)+h_k(s)+\i\mu_k(s)
\]
where $\gamma_k$, $h_k$, and $\mu_k$ are the restrictions of the functions $\gamma$, $h$, and $\mu$ to the interval $J_k=[0,2\pi]$. Note that, by~\eqref{eq:gam}, $\gamma_k(s)=\Im\eta_k(s)$ is known and $h_k(s)=h_k$ where the constant $h_k$ is also known. However, the value of $\mu_k(s)$ will be known only if $s$ is one of the discretization points $s_{k,i}$, $i=1,\ldots,n$, given by~\eqref{eq:sji}. If $s$ is not one of these points, then the value of $\mu_k(s)$ can be approximated using the Nystr{\"o}m interpolation formula (see~\cite{Atk97}). Here, we will approximate $\mu_k(s)$ by finding a trigonometric interpolation polynomial that interpolates the function $\mu_k$ at its known values $\mu_k(s_{k,i})$, $i=1,\ldots,n$. The polynomial can be then used to approximate $\mu_k(s)$ for any $s\in[0,2\pi]$. Finding this polynomial and computing its values is done using the fast Fourier transform (FFT).

\begin{figure}[htb] %
	\centerline{
		\scalebox{0.6}{\includegraphics[trim=0 0 0 0,clip]{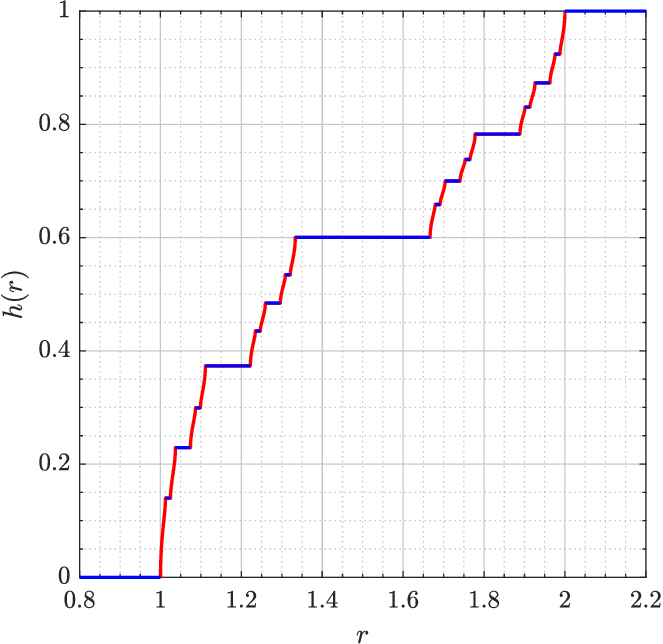}}
		\hfill
		\scalebox{0.6}{\includegraphics[trim=3 0 0 0,clip]{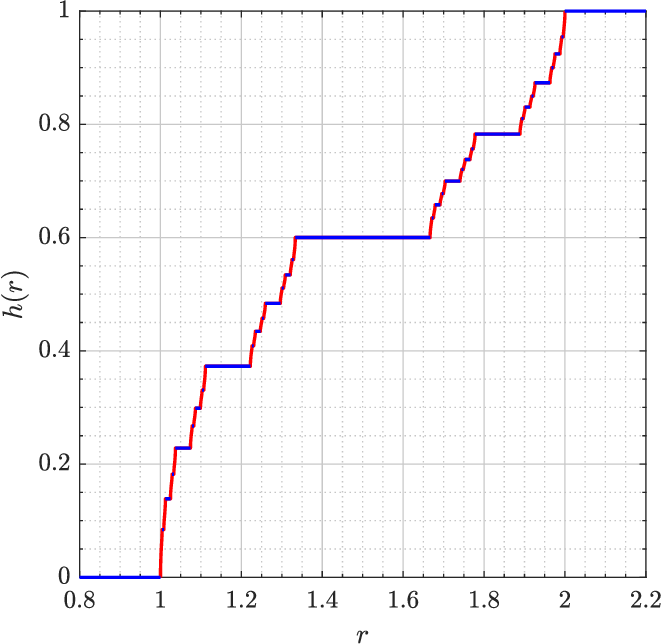}}
	}
	\caption{The $h$-function for the domain $\Omega_{16}$ (left) and the domain $\Omega_{32}$ (right) when $z_0=-3/2$.}
	\label{fig:h-3}
\end{figure}

In our numerical computations below, to compute the non-constant components of the $h$-function, i.e., to compute the values of $h(r)$ when the capture circle intercept with a slit $I_k$ for $k=1,\ldots,m$, we choose $31$ equidistant values of $\xi$ on the upper half of the circle $C_k$. For each of these points $\xi$, we compute the value of $r$ through~\eqref{eq:r-xi} and the values of $h(r)$ through~\eqref{eq:hr3-xi}. 
The values of $h(r)$ for $r$ corresponding to the capture circle passing through the gaps between the slits $I_k$, $k=1,\ldots,m$, are computed as described earlier in Section~\ref{sec:step}.
The graphs of the function $h(r)$ for $m=16$ and $m=32$ are given in Figure~\ref{fig:h-3}.

\subsection{Basepoint $z_0=0$}

\begin{figure}[htb] %
	\centerline{
		\scalebox{0.4}{\includegraphics[trim=0 0 0 0,clip]{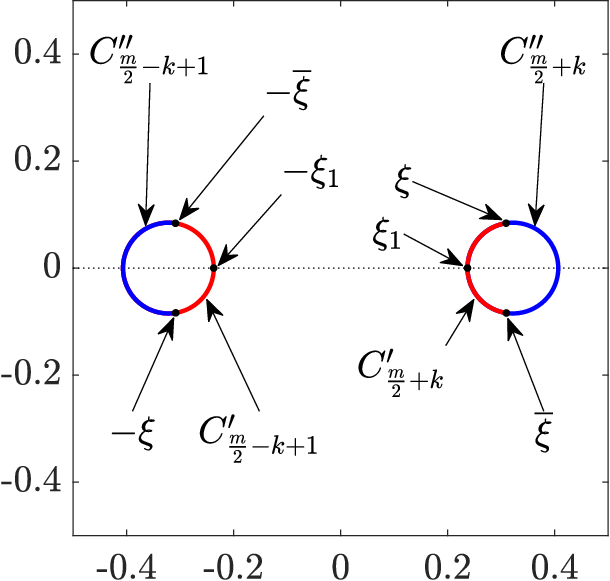}}
	}
	\caption{The arcs $C'_{\frac{m}{2}+k}$, $C''_{\frac{m}{2}+k}$, $C'_{\frac{m}{2}-k+1}$ and $C''_{\frac{m}{2}-k+1}$ for $z_0=0$.}
	\label{fig:h0d}
\end{figure}

For the basepoint $z_0=0$, we assume that the capture circle intercepts the two slits $I_{\frac{m}{2}-k+1}$ and $I_{\frac{m}{2}+k}$ for some $k=1,\ldots,\frac{m}{2}$. Again, due to the symmetry of $G_m$, let $\xi$ and $\overline{\xi}$ be the pair of preimages of the intersection of the capture circle with the slit $I_{\frac{m}{2}+k}$. Also let $\xi_1$ be the intersection of the circle $C_{\frac{m}{2}+k}$ with the positive real axis lying on the arc between $\overline{\xi}$ and $\xi$ which is closest to the point $\zeta_0$. It follows from the symmetry of $G_m$ and $\Omega_m$ that $-\overline{\xi}$ and $-\xi$ are the corresponding pair of preimages of the intersection of the capture circle with the circle $C_{\frac{m}{2}-k+1}$. See Figure~\ref{fig:h0d}.

For $k=1,\ldots,m/2$, the function $U_k(\zeta)$ is the unique solution of the following Dirichlet problem:
\begin{subequations}\label{eq:bdv-U0}
	\begin{align}
		\label{eq:U-Lap}
		\nabla^2 U_k(\zeta) &= 0 \quad \mbox{if }\zeta\in G_m, \\
		\label{eq:U-m}
		U_k(\zeta)&= 1 \quad \mbox{if }\zeta\in C_{\frac{m}{2}-j+1}\cup C_{\frac{m}{2}+j}, \quad j=1,\ldots,k-1, \\
		\label{eq:U-k'}
		U_k(\zeta)&= 1 \quad \mbox{if }\zeta\in C'_{\frac{m}{2}-k+1}\cup C'_{\frac{m}{2}+k},  \\
		\label{eq:U-k''}
		U_k(\zeta)&= 0 \quad \mbox{if }\zeta\in C''_{\frac{m}{2}-k+1}\cup C''_{\frac{m}{2}+k},  \\
		\label{eq:U-p}
		U_k(\zeta)&= 0 \quad \mbox{if }\zeta\in C_{\frac{m}{2}-j+1}\cup C_{\frac{m}{2}+j}, \quad j=k+1,\ldots,m/2. 
	\end{align}
\end{subequations}
Here, $C'_{\frac{m}{2}+k}$ is the arc joining $\overline{\xi}$, $\xi_1$, and $\xi$ on the circle $C_{\frac{m}{2}+k}$ and $C''_{\frac{m}{2}+k}$ is the adjacent arc; $C'_{\frac{m}{2}-k+1}$ is the arc joining $-\overline{\xi}$, $-\xi_1$, and $-\xi$ on the circle $C_{\frac{m}{2}-k+1}$ and $C''_{\frac{m}{2}-k+1}$ is the adjacent arc. The function $U_k(\zeta)$ is assumed to be bounded at infinity.

Note that the boundary conditions on the two circles $C_{\frac{m}{2}-k+1}$ and $C_{\frac{m}{2}+k}$ are again not continuous. In a similar fashion to the problem~\eqref{eq:bdv-U3}, we reduce the current problem~\eqref{eq:bdv-U0} to one with continuous boundary data. 
For $k=1,\ldots,m/2$, the function $\Psi_k(\zeta)$ defined by~\eqref{eq:Psi} is harmonic everywhere in the domain exterior to the $m$ discs, equal to 1 on the arc $C'_{\frac{m}{2}+k}$, and equal to $0$ on the arc $C''_{\frac{m}{2}+k}$.
Similarly, the function 
\begin{equation}\label{eq:Phi_k}
	\Phi_k(\zeta)=\frac{1}{\pi}\Im \log \phi(\psi(\zeta,-\xi,-\xi_1)), 
\end{equation}
is harmonic everywhere in the domain exterior to the $m$ discs, equal to 1 on the arc $C'_{\frac{m}{2}-k+1}$, and equal to $0$ on the arc $C''_{\frac{m}{2}-k+1}$.
Thus, the function 
\begin{equation}\label{eq:Vk-Uk-0}
	V_k(\zeta)=U_k(\zeta)-\Psi_k(\zeta)-\Phi_k(\zeta)
\end{equation}
is a solution of the following Dirichlet problem:
\begin{subequations}\label{eq:bdv-V0}
	\begin{align}
		\label{eq:V-Lap}
		\nabla^2 V_k(\zeta) &= 0 ~~~~~~~~~~~~~~~~~~~~~~~\:\; \mbox{if }\zeta\in G_m, \\
		\label{eq:V-m}
		V_k(\zeta)&= 1-\Psi_k(\zeta)-\Phi_k(\zeta) ~~ \mbox{if }\zeta\in C_{\frac{m}{2}-j+1}\cup C_{\frac{m}{2}+j}, ~ j=1,\ldots,k-1, \\
		\label{eq:V-k'}
		V_k(\zeta)&= -\Phi_k(\zeta) ~~~~~~~~~~~~~~~\:\: \mbox{if }\zeta\in C_{\frac{m}{2}+k},  \\
		\label{eq:V-k''}
		V_k(\zeta)&= -\Psi_k(\zeta) ~~~~~~~~~~~~~~~\:\: \mbox{if }\zeta\in C_{\frac{m}{2}-k+1},  \\
		\label{eq:V-p}
		V_k(\zeta)&= -\Psi_k(\zeta)-\Phi_k(\zeta)  ~~~~~ \mbox{if }\zeta\in C_{\frac{m}{2}-j+1}\cup C_{\frac{m}{2}+j}, ~ j=k+1,\ldots,m/2. 
	\end{align}
\end{subequations}
The function $V_k$ is bounded at infinity.
The boundary data of the new problem~\eqref{eq:bdv-V0} is now continuous on all circles. 

As before, it can be observed also that the image circles under the mappings $\phi(\psi(\zeta,\xi,\xi_1))$ and $\phi(\psi(\zeta,-\xi,-\xi_1))$ are very small, and hence in problem~\eqref{eq:bdv-V0}, the boundary data is almost constant. Hence, in view of (\ref{eq:Psi}) and (\ref{eq:Phi_k}), we may approximate
\begin{equation}
	\Psi_k(\zeta) \approx P_{kj}= \Psi_k(c_j), \quad  \Phi_k(\zeta) \approx Q_{kj}=\Phi_k(c_j), \quad \zeta \in C_j, \quad j=1,\ldots,m,
\end{equation} 
where $c_j$ is the center of the circle $C_j$. Thus, the Dirichlet problem~\eqref{eq:bdv-V0} becomes
\begin{subequations}\label{eq:bdv-V2}
	\begin{align}
		\label{eq:V2-Lap}
		\nabla^2 V_k(\zeta) &= 0 ~~~~~~~~~~~~~~~~~~~~~~~~~~~~~~~~~\:\:\; \mbox{if }\zeta\in G_m, \\
		\label{eq:V2-m}
		V_k(\zeta)&= 1-P_{k,\frac{m}{2}-j+1}-Q_{k,\frac{m}{2}-j+1} ~~ \mbox{if }\zeta\in C_{\frac{m}{2}-j+1}, ~ j=1,\ldots,k-1, \\
		\label{eq:V2+m}
		V_k(\zeta)&= 1-P_{k,\frac{m}{2}+j}-Q_{k,\frac{m}{2}+j} ~~~~~~~~ \mbox{if }\zeta\in C_{\frac{m}{2}+j}, ~ j=1,\ldots,k-1, \\
		\label{eq:V2-k'}
		V_k(\zeta)&= -Q_{k,\frac{m}{2}+k} ~~~~~~~~~~~~~~~~~~~~~~~~ \mbox{if }\zeta\in C_{\frac{m}{2}+k},  \\
		\label{eq:V2-k''}
		V_k(\zeta)&= -P_{k,\frac{m}{2}-k+1} ~~~~~~~~~~~~~~~~~~~~\:\: \mbox{if }\zeta\in C_{\frac{m}{2}-k+1},  \\
		\label{eq:V2-p}
		V_k(\zeta)&= -P_{k,\frac{m}{2}-j+1}-Q_{k,\frac{m}{2}-j+1}  ~~~~\: \mbox{if }\zeta\in C_{\frac{m}{2}-j+1}, ~ j=k+1,\ldots,m/2, \\
		\label{eq:V2+p}
		V_k(\zeta)&= -P_{k,\frac{m}{2}+j}-Q_{k,\frac{m}{2}+j}  ~~~~~~~~~~\; \mbox{if }\zeta\in C_{\frac{m}{2}+j}, ~ j=k+1,\ldots,m/2. 
	\end{align}
\end{subequations}
The function $V_k$ can be then written in terms of the harmonic measures $\{\sigma_j\}_{j=1}^m$ through
\begin{eqnarray}
	\nonumber	V_k(\zeta)&=&
	\sum_{j=1}^{k-1}\left(\sigma_{\frac{m}{2}-j+1}(\zeta)+\sigma_{\frac{m}{2}+j}(\zeta)\right)\\
	\label{eq:Vk-0}	&-&
	\sum_{\begin{subarray}{c} j=1\\j\ne k\end{subarray}}^{m/2} \left([P_{k,\frac{m}{2}-j+1}+Q_{k,\frac{m}{2}-j+1}]\sigma_{\frac{m}{2}-j+1}(\zeta)
	+[P_{k,\frac{m}{2}+j}+Q_{k,\frac{m}{2}+j}]\sigma_{\frac{m}{2}+j}(\zeta)\right)\\
	\nonumber	&-&P_{k,\frac{m}{2}-k+1}\sigma_{\frac{m}{2}-k+1}(\zeta)-Q_{k,\frac{m}{2}+k}\sigma_{\frac{m}{2}+k}(\zeta).
\end{eqnarray}
Then the unique solution $U_k(\zeta)$ to the Dirichlet problem~\eqref{eq:bdv-U0} can be computed through~\eqref{eq:Vk-Uk-0}.

Note that the images of the circles $C_{\frac{m}{2}-k+1}$ and $C_{\frac{m}{2}+k}$ under the conformal mapping $F$ are the slits $I_{\frac{m}{2}-k+1}$ and $I_{\frac{m}{2}+k}$, respectively. Due to the symmetry of $\Omega_m$ and $G_m$ with respect to the real line, we have $F(\xi)=F(\overline{\xi})\in I_{\frac{m}{2}+k}$ and $F(-\xi)=F(-\overline{\xi})\in I_{\frac{m}{2}-k+1}$ (see Figure~\ref{fig:h0d}).
For a given $r$ such that the capture circle intersects the two slits $I_{\frac{m}{2}+k}$ and $I_{\frac{m}{2}-k+1}$ for some $k=1,...,m/2$, there are points $\xi\in C_{\frac{m}{2}+k}$ and $-\xi\in C_{\frac{m}{2}-k+1}$ which depend on $r$ via 
\begin{equation}\label{eq:r0-xi}
	r= F(\xi).
\end{equation}
The values of $F(\xi)$ can be computed as in the case of $z_0=-3/2$. Thus, the value of the $h$-function is given by
\begin{equation}\label{eq:hr0-xi}
	h(r)=U_k(\zeta_0).
\end{equation}

\begin{figure}[htb] %
	\centerline{
		\scalebox{0.6}{\includegraphics[trim=0 0 0 0,clip]{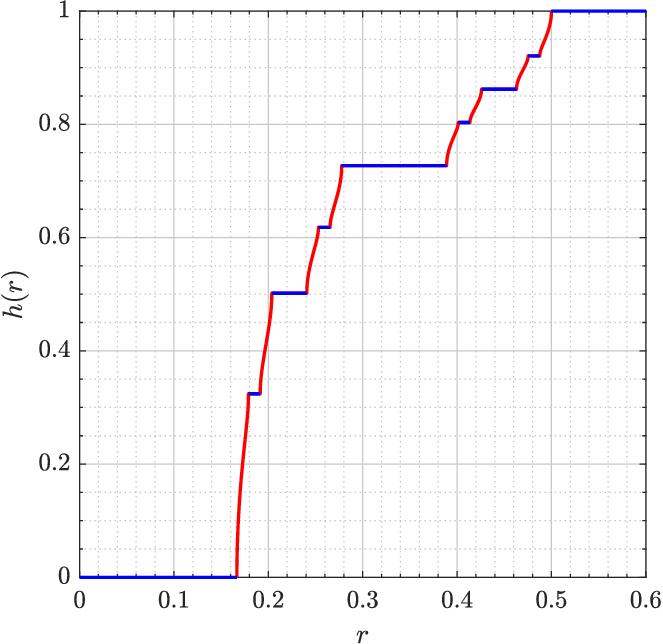}}
		\hfill
		\scalebox{0.6}{\includegraphics[trim=0 0 0 0,clip]{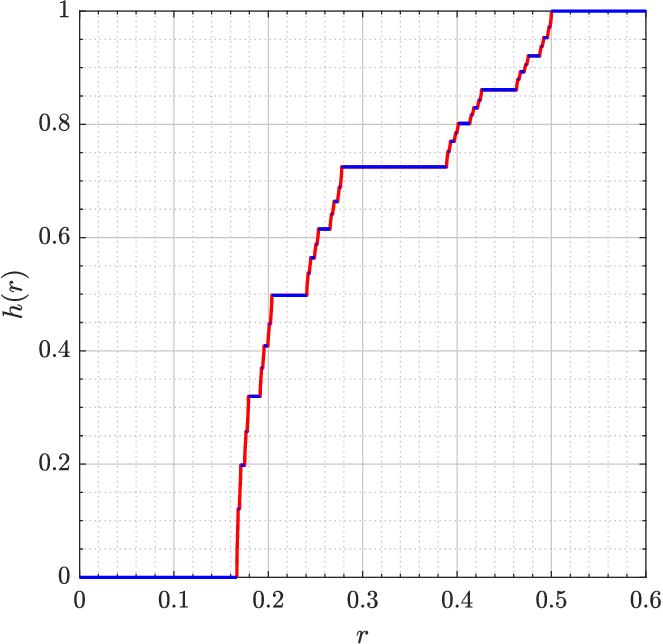}}
	}
	\caption{The $h$-function for the domain $\Omega_{16}$ (left) and the domain $\Omega_{64}$ (right) when $z_0=0$.}
	\label{fig:h-0}
\end{figure}

To compute the values of $h(r)$ when the capture circle intersects the two slits $I_{\frac{m}{2}+k}$ and $I_{\frac{m}{2}-k+1}$ for some $k=1,\ldots,m/2$, as before we take $\xi$ on the circle $C_{\frac{m}{2}+k}$ and then compute the values of $r$ through~\eqref{eq:r0-xi}. 
Again we choose $31$ values of $\xi$ on the upper half of the circle $C_{\frac{m}{2}+k}$. For each of these points $\xi$, we compute the values of $r$ through~\eqref{eq:r0-xi} and then the values of $h(r)$ through~\eqref{eq:hr0-xi}. 
For the values of $r$ corresponding to the capture circle passing through the gaps between the slits $I_k$, the values of $h(r)$ are computed as described in Section~\ref{sec:step}.
The graphs of the function $h(r)$ for $m=16$ and $m=64$ are given in Figure~\ref{fig:h-0}.

The CPU time (in seconds) required for computing the step height of the $h$-function and the $h$-function for $z_0=-3/2$ and $z_0=0$ is presented in Table~\ref{tab:time}. All computations in this paper are performed in {\sc Matlab} R2022a on an MSI desktop with AMD Ryzen 7 5700G, 3801 Mhz, 8 Cores, 16 Logical Processors, and 16 GB RAM. 

\begin{table}[h]
	\caption{The CPU time (sec) required for computing the step heights of the $h$-function (left) and the total time required for computing the $h$-function (right).}
	\label{tab:time}
	\centering
	\begin{tabular}{l|cc|cc}  \hline
		& \multicolumn{2}{c|}{The step height}  & \multicolumn{2}{c}{The $h$-function}  \\ \cline{2-5}
		& $z_0=-3/2$  & $z_0=0$    & $z_0=-3/2$  & $z_0=0$  \\ \hline
		$\Omega_4$       & $0.4120$    & $0.4066$   & $0.4345$    & $0.4363$\\
		$\Omega_8$       & $0.5007$    & $0.5117$   & $0.5186$    & $0.5316$  \\
		$\Omega_{16}$    & $0.8826$    & $0.9272$   & $0.9227$    & $0.9619$  \\
		$\Omega_{256}$   & $24.622$    & $24.715$   & $29.251$    & $29.212$  \\
		$\Omega_{512}$   & $87.541$    & $85.479$   & $104.76$    & $103.08$  \\
		$\Omega_{1024}$  & $330.53$    & $326.17$   & $397.72$    & $392.40$  \\
		\hline
	\end{tabular}
\end{table}

\section{Asymptotic behavior of the $h$-function}\label{sec:asy}

To complement the preceding computations, in this section we analyze numerically some asymptotic features of the $h$-functions as $r\to1^+$ for the basepoint $z_0=-3/2$ and as $r\to(1/6)^+$ for the basepoint $z_0=0$. 

\subsection{The basepoint $z_0=-3/2$}

\begin{figure}[htb] %
	\centerline{
		\scalebox{0.5}{\includegraphics[trim=0 0 0 0,clip]{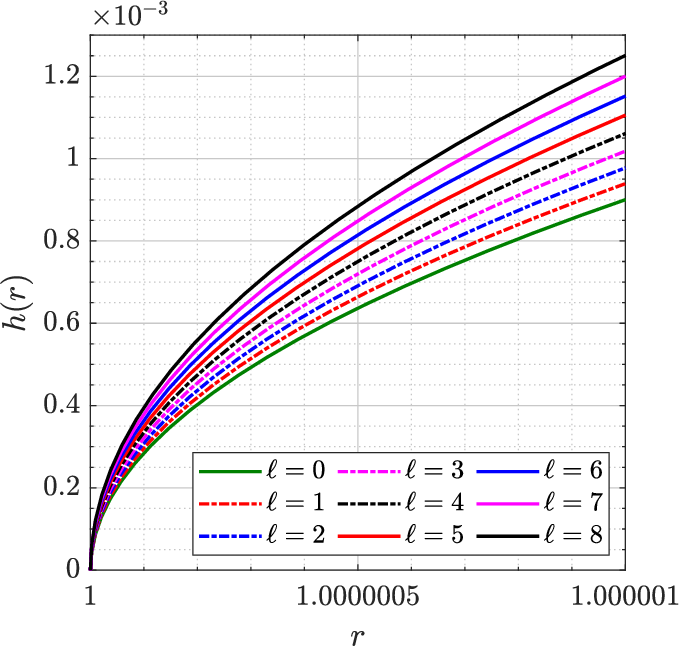}}
		\hfill
		\scalebox{0.5}{\includegraphics[trim=0 0 0 0,clip]{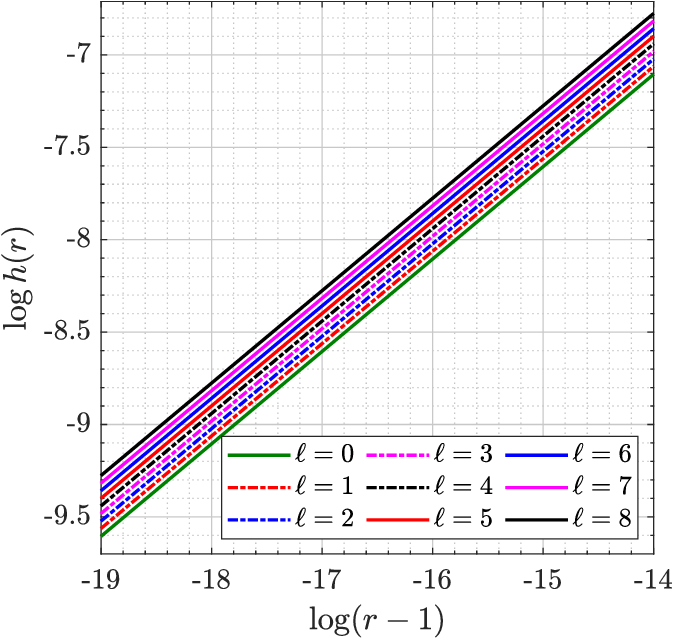}}
	}
	\caption{The graphs of the $h$-function (left) and  $\log(h(r))$ (right) for $\ell=0,1,\ldots,8$ and $r\in(1,1+\varepsilon)$ with $\varepsilon=10^{-6}$ for $z_0=-3/2$.}
	\label{fig:h-asy}
\end{figure}

For this case, we study numerically the behavior of the $h$-function $h(r)$ when $r$ is near to $1$ for several values of $\ell$. 
For $E_0=[-1/2,1/2]$, i.e., $\ell=0$, we have~\cite{gswc}
\begin{equation}\label{eq:h-ell0}
	h(r)=\frac{2}{\pi}\tan^{-1}\left(\sqrt{2}\sqrt{\frac{r-1}{2-r}}\right)
\end{equation}
which implies that
\[
h(r)\sim C_0(r-1)^{\beta_0}, \qquad C_0=\frac{2\sqrt{2}}{\pi}, \quad \beta_0=\frac{1}{2}.
\]
For $\ell\ge1$, up to the best of our knowledge, there are no results about the asymptotic behavior of the $h$-function as $r\to1^+$. Here, we attempt to address this issue numerically.
We choose $20$ values of $r$ in the interval $(1,1+\varepsilon)$ with $\varepsilon=10^{-6}$ (as illustrated in Figure~\ref{fig:3hr}, these $20$ values of $r$ are chosen by choosing $20$ points on the upper-half of the circle $C_1$ such that the images of these points under the conformal mapping $F$ are in  $(1,1+\varepsilon)$). Then, we use the method described in the previous section to compute the values of the $h$-function $h(r)$ for these values of $r$. 
We consider only the values of $\ell=1,\ldots,8$. To consider larger values of $\ell$, it is required to consider smaller values of $\varepsilon$ which is challenging numerically. For $\ell=1,\ldots,8$, the graphs of the functions $h(r)$ on $(1,1+\varepsilon)$ are shown in Figure~\ref{fig:h-asy} (left). The figure shows also the graph of the function $h(r)$ for $\ell=0$ which is given by~\eqref{eq:h-ell0}.
Figure~\ref{fig:h-asy} (right) shows the graphs of the function $\log(h(r))$ as a function of $\log(r-1)$ for $\ell=0,1,\ldots,8$.

Figure~\ref{fig:h-asy} suggests that the functions $h(r)$ for $\ell=1,\ldots,8$ have the same behavior near $r=1$ as for the case $\ell=0$. 
For each $\ell$, $\ell=1,\ldots,8$, we use these $20$ values of $r$, i.e., $r_1,\ldots,r_{20}$, and the approximate values of the $h$-function $h(r)$ at these values to find approximations of real constants $C_\ell$ and $\beta_\ell$ such that
\begin{equation}\label{eq:hrj-C}
	h(r_j)\approx C_\ell (r_j-1)^{\beta_\ell}.
\end{equation}
The approximations of $C_\ell$ and $\beta_\ell$ will be computed using the least square method, i.e., we find the values of $C_\ell$ and $\beta_\ell$ that minimize the least square error 
\begin{equation}\label{eq:LSE}
	\mathcal{E}_\ell = \sum_{j=1}^{20}\left(h(r_j)-C_\ell(r_j-1)^{\beta_\ell}\right)^2.
\end{equation}
By taking the logarithm of both sides of~\eqref{eq:hrj-C}, we obtain
\[
\log(h(r_j))\approx \log(C_\ell)+\beta_\ell\log(r_j-1),
\]
which can be written as 
\[
y_j\approx\log(C_\ell)+\beta_\ell x_j, \quad j=1,\ldots,20,
\]
where $y_j= \log(h(r_j))$ and $x_j=\log(r_j-1)$. The approximate values of the constants $\log(C_\ell)$ and $\beta_\ell$ will be computed by using the {\sc Matlab} function \verb|polyfit| to find the best line fitting for the points $(x_j,y_j)$, $j=1,\ldots,20$. The values of these constant as well as the error $\mathcal{E}_\ell$ in~\eqref{eq:LSE} are given in Table~\ref{tab:Error} for $\ell=1,\ldots,8$.

\begin{table}[h]
	\caption{The values of the constants $C_\ell$ and $\beta_\ell$, as well as the error $\mathcal{E}_\ell$ for  $z_0=-3/2$ (left) and  $z_0=0$ (right).}
	\label{tab:Error}
	\centering
	\begin{tabular}{l|lll|lll}  \hline
		& \multicolumn{3}{c|}{$z_0=-3/2$}   & \multicolumn{3}{c}{$z_0=0$}  \\ \cline{2-7}
		$\ell$   & $C_\ell$      & $\beta_\ell$    & $\mathcal{E}_\ell$     & $C_\ell$  & $\beta_\ell$  & $\mathcal{E}_\ell$ \\ 
		\hline
		$0$      & $0.900316$   & $0.5$         & ---                    & ---         & ---         & ---         \\
		$1$      & $0.939343$   & $0.500000$    & $2.71\times10^{-21}$   & $2.351932$  & $0.500000$  & $9.00\times10^{-18}$  \\
		$2$      & $0.977556$   & $0.500000$    & $1.31\times10^{-20}$   & $2.395871$  & $0.500000$  & $5.42\times10^{-18}$  \\
		$3$      & $1.018398$   & $0.500000$    & $1.89\times10^{-19}$   & $2.466099$  & $0.500000$  & $3.14\times10^{-18}$  \\
		$4$      & $1.061124$   & $0.500000$    & $2.67\times10^{-18}$   & $2.555452$  & $0.500000$  & $4.37\times10^{-19}$  \\
		$5$      & $1.105679$   & $0.500002$    & $1.08\times10^{-17}$   & $2.655781$  & $0.500001$  & $6.70\times10^{-17}$  \\
		$6$      & $1.152042$   & $0.500000$    & $5.67\times10^{-16}$   & $2.763722$  & $0.500004$  & $8.10\times10^{-16}$  \\
		$7$      & $1.200444$   & $0.500001$    & $4.72\times10^{-15}$   & $2.878107$  & $0.500011$  & $1.07\times10^{-14}$  \\
		$8$      & $1.251569$   & $0.500037$    & $2.03\times10^{-14}$   & $2.998958$  & $0.500032$  & $1.31\times10^{-13}$   \\
		\hline
	\end{tabular}
\end{table}

\begin{figure}[htb] %
	\centerline{
		\scalebox{0.5}{\includegraphics[trim=0 0 0 0,clip]{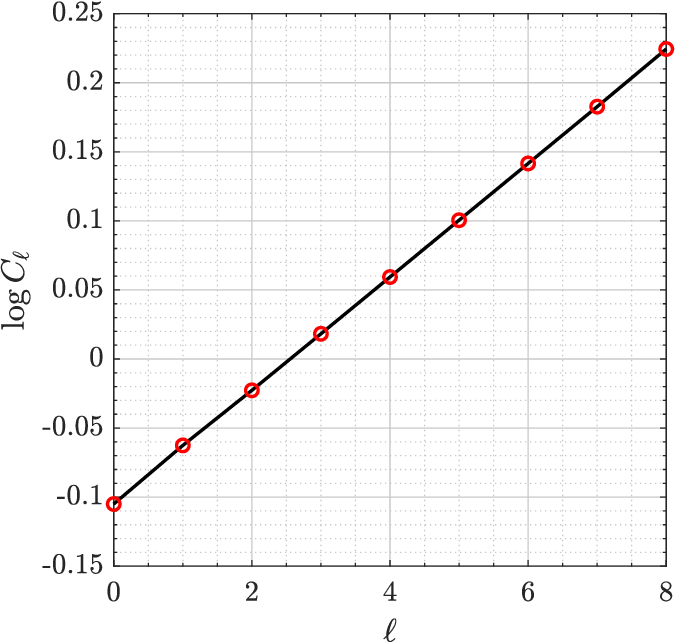}}
		\hfill
		\scalebox{0.5}{\includegraphics[trim=0 0 0 0,clip]{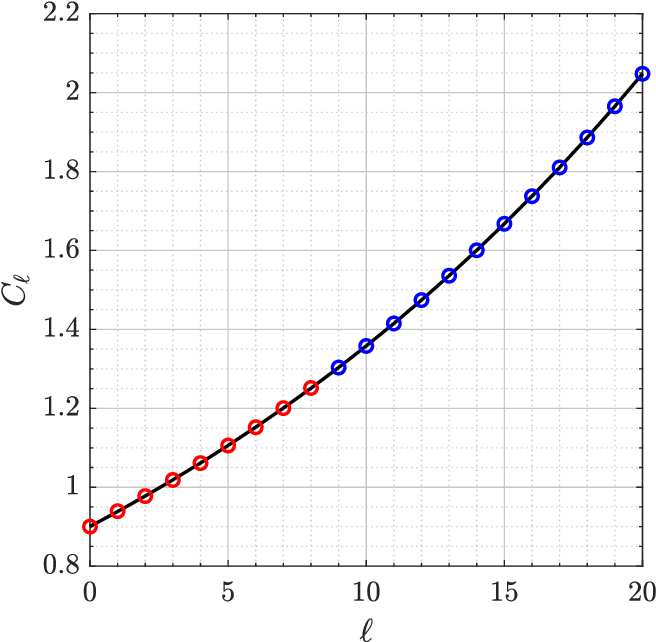}}
	}
	\caption{The values of $\log C_\ell$ (left) and $C_\ell$ (right) for the case $z_0=-3/2$.}
	\label{fig:C08}
\end{figure}

It is clear from Table~\ref{tab:Error} that the values of the constant $C_\ell$ increases as $\ell$ increases. However, it is not clear from the presented values if $C_\ell$ is bounded as $\ell\to\infty$. As a numerical attempt to study the behavior of $C_\ell$ for large values $\ell$, we use extrapolation to find approximate values of $C_\ell$ for $\ell>8$. 
To use the exact value $C_0$ and the computed values $C_1,\ldots,C_8$ to extrapolate the values of $C_\ell$ for $\ell>8$, we note that $\log(C_\ell)$ behave linearly as a function $\ell$ (see Figure~\ref{fig:C08} (left)).   
We again use the {\sc Matlab} function \verb|polyfit| for computing the best line $a\ell+b$ that fits the points $(\ell,\log C_\ell)$ for $\ell=0,1,\ldots,8$. By computing the approximate values of the coefficients, the values of $C_\ell$ can be approximated by
\begin{equation}\label{eq:C-ell}
	C_\ell \approx 0.900613 e^{0.041069\ell}\quad {\rm for}\quad \ell\ge 0.
\end{equation}
The least square error in this formula is
\[
\sum_{\ell=0}^{8} \left(C_\ell - 0.900613 e^{0.041069\ell}\right)^2\approx 1.78\times10^{-6}.
\]

The formula~\eqref{eq:C-ell} suggests that $C_\ell$ is unbounded as $\ell\to\infty$. 
The graph of the $C_\ell$ as a function $\ell$ is shown in Figure~\ref{fig:C08} (right) where the (red) circles are the values of $C_\ell$ presented in Table~\ref{tab:Error}, for $\ell=0,1,\ldots,8$, and the (blue) circles are the extrapolated values, for $\ell=9,\ldots,20$. 

\subsection{The basepoint $z_0=0$}

\begin{figure}[htb] %
	\centerline{
		\scalebox{0.5}{\includegraphics[trim=0 0 0 0,clip]{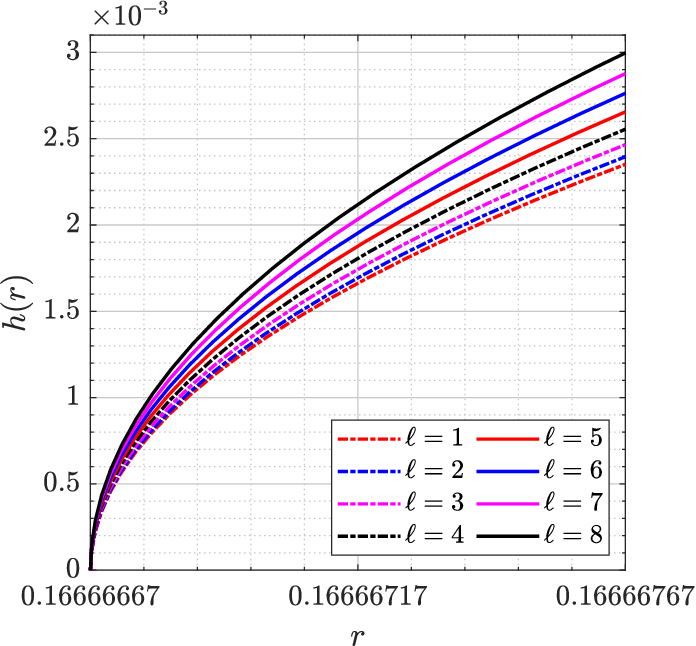}}
		\hfill
		\scalebox{0.5}{\includegraphics[trim=0 0 0 0,clip]{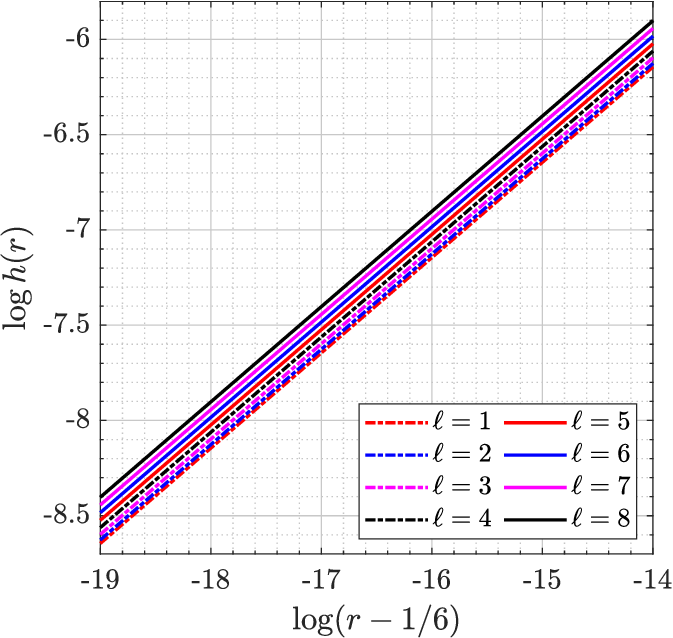}}
	}
	\caption{The graphs of the $h$-function (left) and  $\log(h(r))$ (right) for $\ell=1,\ldots,8$ and $r\in(1/6,1/6+\varepsilon)$ with $\varepsilon=10^{-6}$ for $z_0=0$.}
	\label{fig:h0-asy}
\end{figure}

For this case, it is not possible to consider $\ell=0$ because $z_0=0$ is in the middle of the interval $[-1/2,1/2]$. For $\ell\ge1$, as in the previous case, we choose $20$ values of $r$ in the interval $(1/6,1/6+\varepsilon)$ with $\varepsilon=10^{-6}$ (by choosing $20$ points on the upper-half of the circle $C_1$ such that the images of these points under the conformal mapping $F$ are in  $(1/6,1/6+\varepsilon)$). We then use the method described in Section~\ref{sec:h-fun} to compute the values of the $h$-function $h(r)$ for these values of $r$. 
As in the previous subsection, we consider only the values of $\ell$, $\ell=1,\ldots,8$. The graph of the functions $h(r)$ for these values of $\ell$ on $(1/6,1/6+\varepsilon)$ are shown in Figure~\ref{fig:h0-asy} (left). 
Figure~\ref{fig:h0-asy} (right) shows the graphs of the function $\log(h(r))$ for $\ell=1,\ldots,8$.

Figure~\ref{fig:h0-asy} suggests that the functions $h(r)$ for $\ell=1,\ldots,8$ have the same behavior near $r=1/6$ as for the case $z_0=-3/2$ near $r=1$.
We use the $20$ values of $r$, i.e., $r_1,\ldots,r_{20}\in(1/6,1/6+\varepsilon)$, and the approximate values of the $h$-function $h(r)$ at these values to find approximations of real constants $C_\ell$ and $\beta_\ell$ such that
\begin{equation}\label{eq:hrj-C0}
	h(r_j)\approx C_\ell (r_j-1/6)^{\beta_\ell}
\end{equation}
for each $\ell$, $\ell=1,\ldots,8$.
The approximations of $C_\ell$ and $\beta_\ell$ are computed using the least square method and by using the {\sc Matlab} function \verb|polyfit| as in the case $z_0=-3/2$. The values of these constant as well as the least square error $\mathcal{E}_\ell$ are given in Table~\ref{tab:Error} for $\ell=1,\ldots,8$. By analogy to the case $z_0=-3/2$, it seems from Table~\ref{tab:Error} that the values of $C_\ell$ are also unbounded as $\ell\to\infty$. However, we have been unsuccessful in finding a best fit formula for $C_\ell$ as a function of $\ell$ akin to~\eqref{eq:C-ell} in the case $z_0=-3/2$.

\section{Conclusions}\label{sec:con}

The main objective of this paper was to compute approximations to the $h$-function for the middle-thirds Cantor set $\mathcal{C}$ and to analyze some of its asymptotic features.  
We achieved this by computing $h$-functions for multiply connected slit domains $\Omega_m$ of high connectivity. The computations of the $h$-functions for the slit domains $\Omega_m$ require solving Dirichlet BVPs on these domains. 
We opted to map the slit domains $\Omega_m$ onto conformally equivalent circular domains $G_m$ where the transformed Dirichlet problems are solved and the $h$-functions are then calculated. 
Computing the conformal mapping from $\Omega_m$ onto $G_m$ as well as solving the transformed Dirichlet BVPs in the domain $G_m$ are performed using  the BIE with the Neumann kernel~\eqref{eq:ie-g}. 
This is the same BIE which was used in~\cite{LSN17} to approximate the logarithmic capacity of the middle-thirds Cantor set $\mathcal{C}$.	
As in~\cite{LSN17}, we used the user-friendly {\sc Matlab} function \verb|fbie| from~\cite{Nas-ETNA} to solve the BIE. The FMM lies at the heart of \verb|fbie| and this is what allowed us to consider domains of high connectivity and to perform $h$-function calculations in these domains quickly and accurately.

We have presented results for the step heights of the $h$-function up to connectivity $m=1024$ for the basepoint $z_0=-3/2$ and $m=2048$ for the basepoint $z_0=0$ (see Figures~\ref{fig:hm-3} and~\ref{fig:hm-0}). Our method could certainly be used for larger values of $m$, with still arguably competitive computation times. It is also important to point out that the actual graph of the $h$-function associated with $\mathcal{C}$, regardless of basepoint location, is impossible to compute. However, the main qualitative features of the $h$-functions we have presented for slit domains of high connectivity, such as those in Figures~\ref{fig:hm-3} and~\ref{fig:hm-0}, will differ only very slightly to those exhibited by the actual hypothetical graphs of the $h$-functions associated with $\mathcal{C}$. 
The computed values of the $h$-function have also been used to analyze numerically some asymptotic features of the $h$-functions as $r\to1^+$ for the basepoint $z_0=-3/2$ and as $r\to(1/6)^+$ for the basepoint $z_0=0$. Our work here suggests that the $h$-functions for the multiply connected slit domain $\Omega_m$ ($m>0$) have the same asymptotic profiles as the $h$-function for the simply connected slit domain $\Omega_0$ (see Figure~\ref{fig:h-asy} and Table~\ref{tab:Error}).

Future work will include using the boundary integral equation method to facilitate the calculation of $h$-functions associated with other classes of multiply connected domains and
basepoint locations. Examples of such domains and basepoints have been recently treated in~\cite{ArMa} using the Schottky-Klein prime function. Closely related to the $h$-function is the so-called $g$-function which has received minimal attention to date~\cite{bh89}. It will be interesting to see if our boundary integral equation method can be used to compute $g$-functions for various multiply connected domains and basepoints.

\section*{Acknowledgements}
We thank Lesley Ward and Marie Snipes who suggested to us the topic of this paper. CCG is grateful for the hospitality of the Department of Mathematics, Statistics \& Physics at Qatar University where this work was initiated. CCG also acknowledges support from an Australian Research Council DECRA grant (DE180101098).

\end{document}